\documentclass[11pt,a4paper]{amsart}\usepackage{amssymb,array,bbm}\usepackage{enumerate,graphicx}

\usepackage[bookmarksdepth=3]{hyperref}
\usepackage[utf8]{inputenc}

\newtheorem{Lemma}{Lemma}[section]\newcommand{\bel}{\begin{Lemma}}\newcommand{\eel}{\end{Lemma}}
\newtheorem{Theorem}[Lemma]{Theorem}\newcommand{\bthe}{\begin{Theorem}}\newcommand{\ethe}{\end{Theorem}}
\newtheorem{Proposition}[Lemma]{Proposition}\newcommand{\bprop}{\begin{Proposition}}\newcommand{\eprop}{\end{Proposition}}
\newtheorem{Corollary}[Lemma]{Corollary}\newcommand{\bcor}{\begin{Corollary}}\newcommand{\ecor}{\end{Corollary}}
\newtheorem{Definition}[Lemma]{Definition}\newcommand{\bed}{\begin{Definition}}\newcommand{\eed}{\end{Definition}}
\newtheorem{Definition-Proposition}[Lemma]{Definition-Proposition}

\def\bpr{~\\{\em Proof.\ }}
\newcommand{\epr}{{\hfill\ensuremath\blacksquare}}
\newtheorem{Remark}[Lemma]{Remark}\newcommand{\beR}{\begin{Remark}\rm}\newcommand{\eeR}{\end{Remark}}
\newtheorem{Example}[Lemma]{Example}\newcommand{\bex}{\begin{Example}\rm}\newcommand{\eex}{\end{Example}}
\newtheorem{Problem}[Lemma]{Problem}\newcommand{\bprob}{\begin{Problem}\rm}\newcommand{\eprob}{\end{Problem}}

\newcommand{\beq}{\begin{equation}}\newcommand{\eeq}{\end{equation}}
\newcommand{\bem}{\begin{displaymath}}\newcommand{\eem}{\end{displaymath}}
\newcommand{\beqa}{\begin{eqnarray}}\newcommand{\eeqa}{\end{eqnarray}}
\newcommand{\bee}{\begin{enumerate}}\newcommand{\eee}{\end{enumerate}}
\newcommand{\bei}{\begin{enumerate}[$\bullet$]}\newcommand{\eei}{\end{enumerate}}
\newcommand{\bet}{\begin{tabular}{cccccccc}}\newcommand{\eet}{\end{tabular}}
\newcommand{\bpm}{\begin{pmatrix}}\newcommand{\epm}{\end{pmatrix}}
\newcommand{\bbm}{\begin{bmatrix}}\newcommand{\ebm}{\end{bmatrix}}
\newcommand{\bM}{\begin{matrix}}\newcommand{\eM}{\end{matrix}}
\newcommand{\ber}{\begin{array}{l}}\newcommand{\eer}{\end{array}}

\def\di{\partial}
\def\bl{\langle}\def\br{\rangle}
 
  \newcommand{\isom}[1]{\xrightarrow[\,\smash{\raisebox{1.15ex}{\ensuremath{\scriptstyle\sim}}}\,]{#1}}

 \renewcommand{\stackrel}[2]{\ \lower 0.4ex \hbox{$\mathrel{\mathop{#2}\limits^{\scriptscriptstyle {#1}}}$}\ }

\newcommand{\bin}[2]{\binom{#1}{#2}}

 \newcommand{\quots}[2]{{\footnotesize\left.\raisebox{0.4ex}{$#1$}\! / \!\raisebox{-0.4ex}{$#2$}\right.}}

\def\ca{{\frak a}}
\def\cF{\mathcal{F}}
\def\cO{\mathcal{O}}
\def\cp{{\frak{p}}}
\def\cU{\mathcal{U}}
\def\cm{{\frak m}}

\def\one{{1\hspace{-0.1cm}\rm I}}\def\zero{\mathbb{O}}

\def\iff{if and only if }

\def\k{\mathbbm{k}}\def\N{\mathbb{N}}
\def\P{\mathbb{P}}
\def\R{\mathbb{R}}

\def\Z{\mathbb{Z}}

\def\al{\alpha}\def\de{\delta}
\def\ep{\epsilon}\def\la{\lambda}

\def\tA{\tilde{A}}\def\ta{\tilde{a}}\def\tC{{\widetilde{C}}}
\def\tf{{\tilde{f}}}\def\tJ{{\tilde{J}}}
\def\tl{{\tilde{l}}}\def\tM{\widetilde{M}}\def\tm{\tilde{m}}\def\tn{\tilde{n}}
\def\tQ{\widetilde{Q}}\def\tr{{\tilde{r}}}\def\tS{{\tilde{S}}}
  \def\tU{\tilde{U}}
\def\tZ{{\tilde{Z}}}\def\tz{{\tilde{z}}}

\def\hA{\hat{A}}\def\hR{\hat{R}}
\def\hU{\hat{U}}

\def\empty{\varnothing}
\def\oplusl{\mathop\oplus\limits}\def\liml{\lim\limits}

\def\bA{{\bar{A}}}\def\bB{{\bar{B}}}\def\baf{{\bar f}}\def\bI{{\bar{I}}}\def\bJ{{\overline{J}}}
\def\bbM{{\overline{M}}}\def\bQ{{\overline{Q}}}

\def\bR{{\bar R}}
\def\bU{{\bar{U}}}\def\bV{{\bar{V}}}\def\bx{{\bar x}}
\def\bz{{\bar z}}

\def\smin{\setminus}\def\sset{\subset}\def\sseteq{\subseteq}\def\ssetneq{\subsetneq}

\def\Proj{\operatorname{Proj}}  \def\Spec{\operatorname{Spec}}\def\Tor{\operatorname{Tor}}\def\Coker{\operatorname{Coker}}
\def\mod{mod\text{-}}

\newcommand{\Mat}[1]{\operatorname{Mat}_{m\times n}(#1)}\newcommand{\Mats}[1]{\operatorname{Mat}_{n\times n}(#1)}

\def\uA{{\underline{A}}}
\def\ud{{\underline{d}}}\def\ux{{\underline{x}}}

\def\tf{\tilde{f}}

\title[]{B\MakeLowercase{lock-diagonal reduction of matrices over commutative rings} I.\\
(D\MakeLowercase{ecomposition of modules vs decomposition of their support})}
\author[]{D\MakeLowercase{mitry} K\MakeLowercase{erner and }V\MakeLowercase{ictor} V\MakeLowercase{innikov}}
\address{Department of Mathematics, Ben-Gurion University of the Negev, P.O.B. 653, Be'er Sheva 84105, Israel.}
\email{dmitry.kerner@gmail.com\quad
vinnikov@math.bgu.ac.il}
\thanks{\hspace{-0.3cm}We were   supported by the Israel  Science Foundation (D.K. by grants  1910/18, 1405/22, and V.V.  by grant   2123/17).}
\subjclass[2000]{% Primary 13C14;  Secondary      47A56
 13CXX % Theory of modules and ideals
 13C14 % Cohen-Macaulay modules
 15B33  % Matrices over special rings
 15A21 %Canonical forms, reductions, classification
 15B36 %Matrices of integers
 15A22 % Matrix pencils
 15A54 % Matrices over function rings in one or more variables
 47A56. %Functions whose values are linear operators (operator and matrix valued functions, etc., including analytic and meromorphic ones
}
\date{\today \ filename: \jobname.tex}
\keywords{(Simultaneous) diagonal reduction, (simultaneous) diagonalization, decomposition of modules, 
determinantal representations, determinantal singularities,  matrices over commutative rings, tuples of matrices, matrix version of Noether's $AF+BG$-theorem.
}

\begin{document}
\maketitle  \setcounter{secnumdepth}{4}\setcounter{tocdepth}{1}

\begin{abstract}
Take a rectangular matrix over a commutative ring $A\in \Mat{R}.$ Suppose that the ideal of maximal minors factorizes, $I_m(A)=J_1\cdot J_2\sset R.$
 When is $A$ left-right equivalent to a block-diagonal matrix? (When does the  module/sheaf $\Coker(A)$ decompose as
  $\Coker(A)|_{V(J_1)}\oplus \Coker(A)|_{V(J_2)}$?) If $R$ is not an elementary divisor ring (i.e. not a close relative of a principal ideal ring) one needs additional
   assumptions on $A.$ No necessary and sufficient criterion for such block-diagonal reduction is known.

In this paper we establish the following results.

\bei
\item  The persistence of (in)decomposability under the change of rings. For example:
\bee[-]
\item the passage to Noetherian/local/Henselian/complete rings;
\item  the decomposability of the module $\Coker(A)$ over a graded ring $R$ vs the   decomposability of the sheaf $\Coker(A)$ locally at the points of $\Proj(R)$,
 (this is the matrix decomposability version of M. Noether's $AF+BG$);
\item   the  restriction  to a subscheme in $\Spec(R).$
\eee
\item  The necessary and sufficient condition for decomposability of square matrices in the case: $\det(A)=f_1\cdot f_2$ is not a zero divisor
 and $f_1,f_2$ are co-prime in $R.$
\eei
As an immediate application we give new criteria of simultaneous (block-)diagonal reduction for tuples of matrices over a field,
 i.e., linear determinantal representations.
\end{abstract}
 \tableofcontents

\section{Introduction}
\subsection{}
 Let $R$ be a commutative  unital ring.
Consider the matrices, $A\in \Mat{R}$, $2\le m\le n$,
 up to the left-right equivalence, $A{\sim} U\cdot A\cdot V^{-1}$,  here $U\in GL(m,R)$, $V\in GL(n,R).$
 If $R$ is a principal ideal ring (PIR) then any matrix is equivalent to a diagonal one. This is the well-known Smith normal form.
  More generally, this holds for the elementary divisor rings. These rings were studied in numerous works,
   they are all close to being of Krull dimension one.
   (See, e.g., \cite{Brandal} and  \cite{Karr-Wiegand}.)
 The problem over non-commutative rings was addressed, e.g., in \cite{Ar.Go.O'M.Pa}.
   See also \cite{H.K.K.W.,Wiegand01,Vamos.Wiegand,Zabavsky}.

 For rings of larger Krull dimension, e.g., $\k[[x,y]]$, many matrices are not equivalent to block-diagonal. For example, in the square case,
  $A\in \Mats{R}$,
   the determinant  $\det(A)\in R$  can be irreducible, already this obstructs such a block-diagonal
     reduction (see example \ref{Ex.Decomposability.Square.Matrices.Cases}).

 Our paper grew from the question:
\beq\label{Eq.initial.dec.question}\ber
\text{Suppose the ideal of maximal minors factorizes, $I_m(A)\!=\!J_1\!\cdots\! J_r$, for some ideals $J_i\!\sset\! R.$}\\
\text{How to  ensure the reduction $A\! \sim\! \oplus A_i$, with $A_i\!\in\! \operatorname{Mat}_{m_i\times n_i}(R)$, \ $I_{m_i}(A_i)\!=\! J_i$, \ $\sum m_i\!=\!m$?}
\eer\eeq
\mbox{In fact it is convenient to replace the condition $I_{m_i}(A_i)\!\!=\!\! J_i$ by its weaker version $\sqrt{I_{m_i}(A_i)}\!\!=\!\! \sqrt{J_i}$.}

\medskip

One can ask this question also for the congruence, $A\!\stackrel{congr}{\sim}\!UAU^t$,   and for the conjugation,
  $A\!\stackrel{conj}{\sim}\!UAU^{-1},$ here $U\!\in\! GL(n,R).$ More generally, one needs   decomposability criteria of quiver representations over commutative rings.
 (There is a considerable body of results about containment of quiver representations, see, e.g., \cite{Smalo} and the references therein.)

\medskip

Matrices over rings appear frequently  in pure and applied mathematics, as linear matrix families (or tuples of matrices),
 as matrices of differentiable/analytic functions or power series,
 as the presentations of modules/sheaves/vector bundles and their homomorphisms, as integer matrices (or matrices over rings of integers), and so on.
 The determinantal ideals and determinantal representations  have been intensively studied through decades,  see, e.g.,
  \cite{B.C.V.,Bruns-Vetter,K.M-R.,Miro-Roig}.
The question of (block-)diagonalization is among the most basic, natural, and important.
 Surprisingly it has been  untouched for rings of Krull dimension bigger than one, with the only exception \cite{Laksov}.

\subsection{}
This paper is the first in our study of block-diagonal reduction/block-diagonalization of matrices over commutative unital rings.
 Here we present   results in two directions:
\bei
\item ({\em Obstruction to decomposability and change of base ring}) Take a matrix $A\in \Mat{R}$ with $I_m(A)=J_1\cdot J_2$,
or the module $\Coker(A)\in \mod R$ supported on $V(J_1)\cup V(J_2)\sset \Spec(R).$ We identify the obstruction to the decomposition
 $\Coker(A)\cong \Coker(A)|_{V(J_1)}\oplus \Coker(A)|_{V(J_2)}.$
 This obstruction (an $R$-module) is functorial under the change of base ring.

It is often useful to change the ring.
 For example, one wants to  pass to a Noetherian subring $S\sset R$ that contains the
 entries of $A$,
  or to localize (to check the block-diagonalization on stalks), or to pass to Henselization/completion (i.e., to work with matrices of power series), or
  to take a quotient (i.e., to restrict $A$ to a subscheme in $\Spec(R)$). Under   weak (and natural) assumptions we ensure: $A$  is decomposable \iff
   its image (over $S$) is decomposable.

\item ({\em Decomposability in the square case, $A\in \Mats{R}$}) Assuming $\det(A)\!=\!f_1\!\cdots\! f_r\in R$ (non-invertible),
 an obvious necessary condition to ensure the decomposability as in \eqref{Eq.initial.dec.question}
 is the following bound on the ideal of $(n-1)\!\times\!(n-1)$ minors: $I_{n-1}(A)\!\sseteq\! \sum^r_{j=1} (\prod_{i\neq j}f_i).$
 We prove that this condition is also {\em sufficient} under  rather weak assumptions:
  $\det(A)\in R$ is a regular element and the factors $f_1,\dots,f_r$ are pairwise co-prime.
  Geometrically (in the local regular case) the condition is: the hypersurface germs $V(f_i)\!\sset\! \Spec(R)$ have no common components
   and intersect pairwise properly.

 This solves the decomposition problem for Cohen-Macaulay modules supported on such tuples of hypersurfaces.
In the linear algebra language, we extend the   primary decomposition   for matrices over a field to matrices over a ring.
\eei

Our criteria imply: the decomposability of matrices is controlled by their determinantal ideals $\{I_j(A)\}_j.$
 Recall that $\{I_j(A)\}$ are very naive (rough) invariants of $A$, resp. of the corresponding cokernel module/sheaf $\Coker(A).$
 It was a surprise (for us) that the decomposability question can be settled in quite general case  via these ideals.
 We remark that controlling the ideals  $\{I_j(A)\}$
   is much simpler than controlling the module $\Coker(A).$

Due to   space limitations we give only the first immediate applications to the old problem of simultaneous
  diagonalization of tuples of matrices/linear determinantal representations, decomposition of sheaves on reducible curves, see  \S\ref{Sec.Intro.Structure.Contents}.

\medskip

In the subsequent paper, \cite{Kerner-Block-Diag.II}, we use these decomposability criteria to obtain the  criteria for
 decomposability by  conjugation (decomposability of  representations of groups/algebras), by congruence (decomposability of quadratic/skew-symmetric forms),
  and more generally decomposability of quiver representations (over fields and rings).
   The decomposability criterion for rectangular matrices, also in  \cite{Kerner-Block-Diag.II}, is more delicate. It
   involves controlled cohomology of a determinantal complex of the morphism $A\otimes\quots{R}{I_m(A)}.$

\medskip

\noindent Our interest in block diagonal reduction of matrices over commutative rings of higher Krull dimension originated in our study
 of linear determinantal representations of hypersurfaces \cite{Kerner-Vinnikov.Det.Reps.Global}:
\[
f(x_0,x_1,\ldots,x_n) \!=\! \det(x_0A_0\!+\!x_1A_1\!+\!\cdots\!+\!x_nA_n), \text{ where } A_0, \dots,A_n\!\in\! \Mats{\k} \text{ with } \k \text{ a field}.
\]
This is an old topic in algebraic geometry, see,
e.g., \cite{Dolgachev,Beauville}, which is closely related to matrix factorizations \cite{Eisenbud80,B.H.S.,Backelin-Herzog}.
 We refer to \cite{Kerner-Vinnikov.Det.Reps.Global} for detailed references and to \cite{Vinnikov},\cite{C.S.T.} for
some recent developments and relations.
In that setting, a corollary of our results is that
for $f=f_1\cdot f_2$ with $f_1$, $f_2$ relatively prime, a determinantal representation of $f$ is globally
equivalent to the direct sum of determinantal representations of $f_1$ and of $f_2$ \iff
 it decomposes locally at every point of intersection of the corresponding hypersurfaces.
In the case of curves ($n=2$) and assuming that the two curves $f_1=0$ and $f_2=0$ have no common tangents
at every point of intersection, a determinantal representation decomposes \iff it is maximally generated
(the dimension of its kernel equals the multiplicity of the point on the curve $f=0$,
also called Ulrich maximal \cite{Ulrich84}) at every point
of intersection,
see Corollary \ref{Thm.Cor.4.7}
and \S\ref{Sec.Decomposition.Det.Reps.Curves}.

We mention in this context also the recent results \cite{KlVo17,HeKlVo18} on equivalence and decomposition
of matrices of linear forms, when these forms are viewed as polynomials in free noncommuting variables
and are evaluated on square matrices of all sizes over the ground field.

\subsection{The structure/contents of the paper}\label{Sec.Intro.Structure.Contents}
\bee[\S 1\!]\setcounter{enumi}{1}
\item contains the relevant background.
 
In \S\ref{Sec.Background.Notations.Conventions} we set the notation  and define the (stable) decomposability of matrices.

Then we recall the (trivial) decomposability of modules with non-connected support and the openness of decomposability locus.

  \S\ref{Sec.Background.Coprime.elements.Coregular.ideals} gives some homological versions of
 ``co-regularity" of ideals $J_1\cap J_2=J_1\cdot J_2$, e.g., via $Tor_1(\quots{R}{J_1},\quots{R}{J_2})$ and $H^0_{J_1+J_2}(M).$

\item treats the obstruction to decomposability and the reduction to ``convenient" rings.
\bee[\!\!\!\!\S 3.1.\!]
\item We prove: $A$ is decomposable over $R$ \iff $A$ is decomposable over certain Noetherian subring $S\sset R$ that contains the entries of $A.$
 Moreover, (in the square case) for $A\in \Mats{R}$ the assumptions ``$\det(A)=\prod f_i$, with $\{f_i\}$ co-prime, and $I_{n-1}(A)\sseteq \sum^r_{j=1} (\prod_{i\neq j}f_i)$"
  hold over $R$ \iff they hold over $S$, assuming a particular chain stabilization condition, much weaker than Noetherianity.

  \item When the support of a module is reducible, $Supp(\Coker(A))=V(J_1)\cup V(J_2)\sset \Spec(R)$, it is natural to separate the components,
 $V(J_1)\amalg V(J_2)\stackrel{\pi}{\to}V(J_1)\cup V(J_2).$ Then one compares $\Coker(A)$ to (roughly) $\pi_*\pi^*\Coker(A).$ The module $\pi_*\pi^*\Coker(A)$
  is obviously decomposable. Comparing the two modules one gets the obstruction to the decomposability,
  \[
  0\to H^0_{J_1+J_2}(\Coker(A))\to \Coker(A)\to \pi_*\pi^*\Coker(A)\to Q\to0.
  \]

\item  We prove: this obstruction $Q$ is functorial under  those base changes that are compatible with  certain local cohomology objects.
\item
 As an application we get:
 \bei
 \item The matrix $A$ is decomposable over $R$ \iff $A$ is decomposable
  over all the localizations, $R_\cm.$
  \item For \!a local ring,
 $A$ is \!decomposable over $R$ \iff it \!is \!decomposable over the Henselization  $R^h.$
  \item One can   pass to the completion, $ \hR^{(\cm)}$, assuming the completion functor $\otimes \hR^{(\cm)}$ is exact.
\eei

\item For the rings of Analysis, e.g., $C^\infty(\cU)$ or  $C^\infty(\R^n,o)$, the functor $\otimes \hR^{(\cm)}$ is neither exact nor faithful.
 Thus we give a separate argument for the reduction  to
  complete rings.

\item In many cases one needs a non-flat base change $R\to S$, e.g., for $S=\quots{R}{\ca}.$ In this case we give simple conditions to ensure:
 $A$ is decomposable over $R$ \iff  $A\otimes S$ is decomposable.

\item Suppose the ring is  graded,  $R=\oplus_{d\in \N}R_d$, and
  the  matrix is  graded, see \S\ref{Sec.Background.Notations.Conventions}.v.
 Denote by $\cm$ the maximal homogeneous ideal in $R.$
  We prove: $A$ is $GL(m,R)\times GL(n,R)$-decomposable \iff the $\cm$-localization, $A_\cm,$ is
 $GL(m,R_\cm)\times GL(n,R_\cm)$-decomposable.

 \item Suppose $R,A$ are graded. Using the standard correspondence of graded $R$-modules to coherent sheaves, $mod_{gr}$-$R\rightsquigarrow Coh(Proj(R))$, we
  can consider   $\Coker(A)$ as the sheaf of modules on the projective scheme $\Proj(R).$
 The decomposition of $A$ over $R$    is reduced to the
 local decompositions of the stalks of the sheaf $\Coker(A)$  at  the   points of the
   subscheme $\P V(J_1)\cap \P V(J_2)\sset \Proj(R).$ %, proposition \ref{Thm.Decompos.Proj(R).vs.Local}.
    This is the matrix decomposability version of the fundamental $AF+BG$-theorem  by M. Noether.

 Such a reduction is perhaps unexpected, as the scheme $\P V(J_1)\cap \P V(J_2)$ can be non-connected, e.g., a finite set of closed points.
 (Alternatively, one could expect monodromy-type effects.)

This M. Noether-type result gives the traditional reduction in dimension:
 ``A question over a graded ring $R$ is reduced to many questions over local rings of dimension $dim(R)-1$".

\eee
The results of \S3 are used both in \S4 and in \cite{Kerner-Block-Diag.II}.

\item gives the  decomposition criterion for square matrices, theorem \ref{Thm.Decomposition.Square.Matrices}.
 This is a complete solution of the decomposition problem for  matrices with $\det(A)=f_1\cdots f_r$, here $\{f_i\}$ are regular and pairwise co-prime.

 \bee[\!\!\!\!\!\S4.1.\!]
\item (a preparation) We establish the block-diagonalization (by conjugation) of projectors and ``almost projectors" over  local/Henselian rings.
  (By \S\ref{Sec.Reduction.to.local.complete} we can always assume $R$ is local-Henselian, though not necessarily Noetherian.)

\item Theorem \ref{Thm.Decomposition.Square.Matrices} is proved by creating such almost projectors from $A$ and its adjugate $Adj(A).$
  The proof is completely down-to-earth, using only the $R$-linear algebra.

In \cite{Kerner-Block-Diag.II} we give a shorter proof of theorem \ref{Thm.Decomposition.Square.Matrices}, via a certain determinantal complex. But the current down-to-earth proof has its advantages.
 It is easily adapted to decomposability by the congruence, $A\to UAU^t$, and is useful when extending  theorem \ref{Thm.Decomposition.Square.Matrices} to the case of non-commutative  rings.

\item
The first examples. As a trivial application of theorem \ref{Thm.Decomposition.Square.Matrices} we derive the ``first half" of the Smith
 normal form over PID's. Then we give decomposability criteria for determinantal representations of maximal corank (Ulrich-maximal modules).

\item Discussion of the geometry of the assumptions in theorem \ref{Thm.Decomposition.Square.Matrices}.

\item Decomposability criteria for (non-linear) determinantal representations and for torsion-free sheaves on reducible plane curves.

\item We get criteria for simultaneous diagonal reduction of tuples of matrices over a field
 (i.e., linear determinantal representations).
 \bei
 \item
 The case of pairs of matrices is elementary and classically known.
\item The general case is reduced to the case of triples of matrices (preserving the size!), via \S\ref{Sec.Reduction.Passage.to.Quotient}.
 \item Such triples are determinantal representations of projective plane curves. For triples admitting the (full) diagonal reduction the determinantal curve is a line arrangement. The triple of matrices admits the diagonal reduction \iff   ``the total defect" of
  kernel-dimensions  (properly counted) has the maximal possible value.
\eei
No criteria of such type could be imagined with the previous (classical) methods.
  \eee
\eee

\subsection{Acknowledgements}
  The first results of this paper appeared (long ago) as an offshoot of our work on determinantal representations, \cite{Kerner-Vinnikov.Det.Reps.Global}.
 Many of the results were obtained  during the postdoctoral stay of D.K. in  the University of Toronto, 2010-2012, \cite{Kerner-Vinnikov.Decompos.Det.Reps}.

We are thankful to G. Belitski, A. F. Boix, to the late R.O. Buchweitz, to I.
Burban, G.M. Greuel, M. Leyenson, I. Tyomkin, M. Zach for important discussions through all these years.

 Special thanks are to the referee for very helpful remarks.

\section{Preparations}\label{Sec.Background}
\subsection{Notation  and conventions}\label{Sec.Background.Notations.Conventions}
 Unless stated otherwise, $R$ is a commutative unital ring, not necessarily Noetherian, while
$\k$ denotes a field, of any characteristic.
 When writing ``$(R,\cm)$ is a local ring", we do not assume $R$ is Noetherian.  The ideals are not assumed finitely generated.

\bee[\!\!\!\!\bf i.\!]
\item
Denote the (square) unit matrix by $\one$, the zero matrix (possibly non-square) by $\zero.$
 Let $A\in \Mat{R}$, we always assume $2\le m\le n.$

\item
Sometimes we change the ring, $R\to S.$ Then we take the images,   $R\ni f\to \baf\in S$ and
$\Mat{R}\!\ni\! A\to\bA\!\in\! \Mat{S}.$ For an ideal $I\sset R$ denote by $\bI\!\sset\! S$ the ideal generated by the image of $I.$
 (No confusion with the integral closure, as we do not use it in this paper.)

For $M=\Coker(A)$ we denote $\bar{M}=\Coker(\bA).$

Applying $S\otimes-$ to a submodule $M\sset N$ we denote by $S\cdot M$ the image of $S\otimes M\to S\otimes N.$

\item
The  determinantal ideal of a matrix,  $I_j(A)$,  is generated by all the $j\times j$ minors of $A\in \Mat{R}$, see
\cite[\S20]{Eisenbud-book}.
 The determinantal ideals form a decreasing chain,
 \beq
 R=I_0(A)\supsetneq I_1(A)\supseteq\cdots\supseteq I_m(A)\supseteq I_{m+1}(A)=0.
 \eeq
 This chain is invariant under the $GL(m, R)\times GL(n, R)$ action, $I_j(A)=I_j(UAV^{-1}).$

Consider $A$ as a presentation matrix of its cokernel, $R^n\stackrel{A}{\to}R^m\to \Coker(A)\to0.$ Usually we denote $M:=\Coker(A).$
 Then the determinantal ideals coincide with the Fitting ideals: $I_j(A)=Fitt_{m-j}(M).$
 In particular, as $Fitt_0(M)\cdot M=0$, the image of $A$ satisfies: $Im(A)\supseteq I_m(A)\cdot R^m.$

In fact a stronger property holds: $I_{m-1}(A)\cdot Im(A)\supseteq I_m(A)\cdot R^m.$

The determinantal ideals are functorial, for $R\to S$ one has $I_j(\bA)=I_j(A)\cdot S.$
\item
The adjugate of a square matrix $A\in \Mats{R}$ is the matrix of cofactors, $Adj(A)\in\Mats{R}.$
 Usually in this paper  $\det(A)\in R$ is not a zero divisor, then $Adj(A)$ is determined  by the condition $Adj(A)\cdot A=\det(A)\cdot\one=A\cdot Adj(A).$
 The entries of $Adj(A)$ generate the ideal $I_{n-1}(A)$ and one has $\det(Adj(A))=det(A)^{n-1}.$

\item We use two kinds of matrix equivalence:
\bei
\item the   left-right equivalence: $A\sim B$ if $A=UBV^{-1}$ for some $U\in GL(m,R)$, $V\in GL(n,R)$;
\item the {\em stable} left-right equivalence:  $A\stackrel{s}{\sim}B$ if $U(A\oplus \one_{r\times r}) V^{-1}= B\oplus\one_{\tr\times \tr} $ for some $r,\tr\in \Z_{\ge0}$ and
 $U\in GL(m+r,R)$, $V\in GL(n+r,R).$ (Here $A,B$ can be of different sizes.)
\eei
The stable equivalence often implies the ordinary equivalence, e.g., this happens in the following cases.
\bei
\item When $R$ is a local ring.
\item  When $R$ is graded, $R=\oplus_{d\in \N}R_d$, with $R_0$ local, Noetherian, and the matrices $A,B$ are graded,
 i.e.,  their entries are homogeneous and the degrees satisfy $deg(a_{ij})+deg(a_{kl})=deg(a_{il})+deg(a_{kj})$ for all $i,j,k,l.$
   This latter condition can be stated as $\{deg(a_{ij})=d_i+d_j\}_{i,j}$ for some integers $d_\bullet.$
\item  When $R$ is a principal ideal ring.
\eei
 This follows by the uniqueness of projective resolution of the $R$-module $\Coker(A)$, \cite[\S20]{Eisenbud-book}.

Over   a local ring, $(R,\cm),$ one can pass to the minimal
resolution, therefore $A\sim \one\oplus \tA$, with $\tA\in
\operatorname{Mat}_{\tm\times \tn}(\cm).$

\item (The decomposability)
\bed  Let $A\in \Mat{R}$ with the factorized ideal of maximal minors $I_m(A)=J_1\cdots J_r\sset R$, $2\le m\le n.$
\bee[1.]
\item   $A$ is called $(J_1,\dots, J_r)$-decomposable if  $A \sim  \oplus^r_{i=1} A_i$, with $\sqrt{I_{m_i}(A_i)}\!\!=\!\!\sqrt{J_i}\sset R.$
\item
$A$ is called stably-$(J_1,\dots, J_r)$-decomposable if  $A\stackrel{s}{\sim} (\oplus^r_{i=1} A_i)$,
 with $\sqrt{I_{m_i}(A_i)}\!\!=\!\!\sqrt{J_i}\sset R.$
\item
A module $M$ with the factorized  zeroth Fitting ideal, $Fitt_0(M)=J_1\cdots J_r  \sset R,$  is called $(J_1,\dots, J_r)$-decomposable if
  $M= \oplus M_i$, with $\sqrt{Fitt_0(M_i)}= \sqrt{J_i}.$
\eee\eed
\noindent We use here the radicals of   ideals, rather than e.g. $ I_{m_i}(A_i) \!\!=\!\!J_i,$ because the ideals $J_1,J_2$  are not determined uniquely by
 the condition $I_m(A)\!=\!J_1\!\cdot\! J_2,$ even with the additional assumption $J_1\!\cap\! J_2\!=\!J_1\!\cdot\! J_2.$
  (Only the minimal associated primaries of $J_1,J_2$ are uniquely defined.)

 \medskip

 Recall the standard fact:
\beq\label{Eq.Decompos.of.A.vs.Coker(A)}
\text{$A$ is stably-$(J_1,\dots, J_r)$-decomposable \iff $\Coker(A)$ is $(J_1,\dots, J_r)$-decomposable.}
\eeq
This follows directly by Fitting's lemma for the resolutions of the modules $\Coker(A)$, $\oplus \Coker(A_i)$, see \cite[\S A.3]{Eisenbud-book}

\item(Restriction onto the support)
Take a local ring $(R,\cm)$, two matrices  $A,B\in\Mat{\cm}$ with $2\le m\le n$,  and the quotient homomorphism $ R \to \quots{ R}{I_m(A)}.$
 Then $A\sim B$  \iff   $\bA\sim \bB\in \Mat{\quots{ R}{I_m(A)}}.$
\bpr The part $\Rrightarrow$ is trivial, we prove the part $\Lleftarrow.$

   Suppose $\bA\sim \bB$ then   $A=UBV^{-1}+Q$, for some $U\in GL(m, R)$, $V\in GL(n, R)$ and $Q\in \Mat{I_m(A)}.$
 By \S\ref{Sec.Background.Notations.Conventions}.iii we have: $Q=A\cdot Q'$, for a matrix $Q'\in \Mats{\cm}.$
  Therefore $A=UBV^{-1}\cdot (\one-Q')^{-1}.$
\epr

\medskip

\noindent
The locality  is needed e.g., because of the trivial example: $R=\k[x]$, $A=x\one$, $B=x(x-1)\one.$
\\The assumption $2\le m$ is needed because of the trivial example: $B=\zero\in \rm{Mat}_{1\times n}(R)\ni A\neq \zero.$

\medskip

\noindent Therefore the decomposability of $A$ is testable  by restricting to $V(I_m(A)) \sset  \Spec(R)$:
\beq
\hspace{1cm}A\stackrel{s}{\sim}A_1\oplus A_2 \hspace{2cm}\text{ \iff }\hspace{2cm}\bA\stackrel{s}{\sim}\bA_1 \oplus \bA_2\in \Mat{\quots{ R}{I_m(A)}}.
\eeq
\eee

\subsection{Decomposability of modules with non-connected support}\label{Sec.Background.Decomposability.for.non.Connected.Support}
\bel
Let $A\in \Mat{R}$ and assume $I_m(A)=J_1\cdot J_2$, for some ideals $J_1,J_2\ssetneq R$ satisfying $J_1+J_2=R.$
  Then   $A$ is stably $(J_1,J_2)$-decomposable.
\eel
\noindent
Geometrically, if $Supp(M)\!=\!\!V(J_1)\!\amalg\! V(J_2)\!\sset\! \Spec(R)$ then $M$ decomposes, $M\!\cong\!  M|_{V(J_1)}\!\oplus\! M|_{V(J_2)}.$
\bpr
We have: $\Coker(A)=J_1\cdot \Coker(A)+ J_2\cdot \Coker(A).$
This is a direct sum:
\beq
J_1\cdot \Coker(A)\cap J_2\cdot \Coker(A)=(J_1+J_2)\cdot\big[J_1\cdot \Coker(A)\cap  J_2\cdot \Coker(A)\big]=0.
\eeq
 (Note that $J_1J_2\cdot \Coker(A)=0$, by \S\ref{Sec.Background.Notations.Conventions}.iii.)
 \quad \quad Now invoke \eqref{Eq.Decompos.of.A.vs.Coker(A)}.
\epr

\medskip

This statement is for stable decomposability, the ordinary one does not hold, see remark \ref{Rem.3.13}.

\subsection{Openness of the decomposability locus}
\bel
Let $I_m(A)=J_1\cdot J_2$ and suppose  the localization $A_\cm\in \Mat{R_\cm}$ at a maximal ideal $\cm\supseteq J_1+J_2$ is
 $(J_1)_\cm,  (J_2)_\cm$-decomposable. Then there exists an open neighborhood $V(\cm)\in \cU\sseteq \Spec(R)$ such that $A|_\cU$ is decomposable.
\eel
  Namely, there exists a (non-nilpotent) element $g\in R\smin \cm$ such that for the base change $R\to R[\frac{1}{g}]$ the matrix
  $\bA\in Mat_{m\times n}(  R[\frac{1}{g}])$ is $(\bJ_1,\bJ_2)$-decomposable.
  \bpr
  Suppose $U_\cm\cdot A_\cm\cdot V_\cm=A_{1,\cm}\oplus A_{2,\cm}.$ Present the entries of $U_\cm,V_\cm$ as fractions. Let $g\in R\smin \cm$ be the product
   of all the denominators of $U_\cm,V_\cm.$ Take $\cU:=Spec(R)\smin V(g)=Spec(R[\frac{1}{g}]).$ Then
    $\bU\bA\bV=\bA_1\oplus\bA_2\in \Mat{R[\frac{1}{g}]}.$
  \epr

\subsection{Co-prime elements and co-regular ideals}\label{Sec.Background.Coprime.elements.Coregular.ideals}
Two regular elements $f_1,f_2\in R$ (i.e.,  neither invertible, nor zero divisors) are called co-prime if
$(f_1)\cap (f_2)=(f_1\cdot f_2)\sset R.$   This can be stated also as: both $f_1,f_2$ and $f_2,f_1$ are  regular sequences in $R.$
 This is equivalent to: for  every presentation $\{f_i=g_i\cdot h\}_i$ the element $h\in R$ is invertible.

\bel\label{Thm.Coprime.Elements}
\bee[\bf 1.\!]
\item The regular elements $\{f_i\}\sset R$  are
pairwise co-prime   \iff the pairs $\{f_i,\prod_{j\neq i}f_j\}$ are co-prime
for all $i.$

\item Assume $ f_1,\dots, f_r\in R$ are regular and pairwise co-prime.
 Then  $
\sum^r_{j=1} (\prod_{i\neq j} f_i)=\cap^r_{j=1} \big((f_j)+(\prod_{i\neq j} f_i)\big)\sset R.$
\eee
\eel
\bpr Denote $h=\prod f_j,$ thus $\prod_{j\neq i}f_j=\frac{h}{f_i} $ and $\prod_{j\neq i,k}f_j=\frac{h}{f_if_k}.$
\bee[\bf 1.]
\item  {\bf The part $\Rrightarrow$.} If $af_i\in (\frac{h}{f_i})$ then $a=\ta\cdot f_k$, for some $k\neq i$. Therefore
 $\ta f_i\in (\frac{h}{f_if_k})$. And so on.
\\
{\bf The part $\Lleftarrow$.} If  $af_i\!\in\! (f_k)$ then $(\frac{h}{f_if_k})
af_i\!\in\! (\frac{h}{f_i })$. Therefore
 $(\frac{h}{f_if_k}) af_i\!\in\! (h)$, and hence $a\!\in\! (f_k)$.

\item
  The inclusion $\sseteq$ is obvious. For the inclusion $\supseteq$ we denote $h=\prod f_j$
 and  take an element
 $a_k f_k+b\frac{h}{f_k}\in   \cap^r_{j=1} \big((f_j)+(\frac{h}{f_j})\big)$.
 It is enough to prove: $a_kf_k\in \sum^r_{j=1} (\frac{h}{f_j})$.
 By the assumption   $a_k f_k= a_i f_i  +b_i \frac{h}{f_i}  $, thus (as $f_k,f_i$ are co-prime) $a_k= \ta_i f_i +b_i  \frac{h}{f_if_k}  $.
  Thus it is enough to prove: $ \ta_i f_if_k\in   \sum^r_{j=1} (\frac{h}{f_j})$. And so on.
\epr
\eee

\medskip

We often impose this ``co-regularity" condition on general ideals: $J_1\cap J_2\! =\! J_1\cdot J_2\sset R$.
 Geometrically, for $(R,\cm)$ local Noetherian, this implies: the subschemes $V(J_1),V(J_2)\!\sset\! \Spec(R)$ intersect properly and contain no embedded
  components supported on $V(J_1\!+\!J_2)$.
\bel\label{Thm.coprime.ideals}
\bee[\bf 1.]
\item $J_1\cap J_2=J_1\cdot J_2$ \ \iff \ $\Tor_1^R(\quots{R}{J_1},\quots{R}{J_2})=0$.
\item Let $A\in \Mat{R}$.
Suppose $I_m(A)\!=\!J_1\!\cdot\! J_2\!=\!J_1\!\cap\! J_2\!\sset\! R$ and take the $R$-module $M\!:=\!\Coker(A)$.
\bee[\bf i.]
\item  Then $H^0_{J_1+J_2}(M)\supseteq J_1  M \cap  J_2  M$.
\item If  $ J_1  M \cap  J_2  M =0$ then $(J_iM):_M J_j=0:_MJ_j$ for $i\neq j$.
\item (Square-matrix case) Suppose $\det(A)=f_1\cdot f_2$, regular and co-prime. Then $H^0_{(f_1,f_2)}(M)=0$.
\eee
\eee
\eel
\bpr (These facts are standard, we recall the proof.)
\bee[\bf 1.]
\item   Apply $\otimes \quots{R}{J_2}$ to the exact sequence $0\to J_1\to R\to\quots{R}{J_1}\to0$ to get:
\beq\label{Eq.Coprime.ideals.vs.Exact.sequence}
0=\Tor_1^R(R,\quots{R}{J_2})\to \Tor_1^R(\quots{R}{J_1},\quots{R}{J_2}) \to J_1\otimes \quots{R}{J_2}\stackrel{\phi}{\to}  \quots{R}{J_2}\to\quots{R}{J_1+J_2}\to0.
\eeq
Therefore $\Tor_1^R(\quots{R}{J_1},\quots{R}{J_2})=0$ \iff the map $\phi$  is injective \iff  $J_1\cap J_2=J_1\cdot J_2$.
\item
\bee[\!\!\!\!\bf i.]
\item We have $(J_1+J_2)\cdot (J_1  M \cap  J_2  M )=J_1\cdot J_2\cdot M=0$.
\item The inclusion $\supseteq$ is obvious.  For the part $\sseteq$ we note: $J_j\cdot(J_iM:_M J_j)\sseteq J_iM\cap J_jM=0$.
\item It is enough to prove: $0:_M(f_1,f_2)=0$. Take an element of $0:_M(f_1,f_2) $ and its representative $\xi\in R^n.$
  Then $(f_1,f_2)\cdot \xi\sset Im(A)$. Therefore  $(f_1,f_2)\cdot Adj(A)\cdot\xi\sset f_1f_2\cdot R^n$. As each of $f_i$ is regular, we get:
  $Adj(A)\cdot\xi\in (f_1\cdot R^n)\cap (f_2\cdot R^n)=f_1f_2\cdot R^n$. Therefore $\xi\in Im(A)$, i.e.,  $\xi$ represents $0\in M$.
\epr
\eee
\eee

\section{The obstruction to decomposability and reduction to ``convenient" rings}\label{Sec.Reduction}
 When establishing the decomposability conditions it is often useful to change the base ring, e.g., to pass to a subring $S\sset R$
  (that contains the entries of $A$), or to extend/to take the quotient, $R\to S$.
  To trace the behaviour of $A$ and $M:=\Coker(A)$ under this change we establish results of two types:
  \bee[\bf i.]
  \item
  $A$ is decomposable over $R$ \iff $\bA$ is decomposable over $S$.

\item
 Consider the following conditions (they are needed for the decomposability, see \S\ref{Sec.Decomposability}):
\beq\label{Eq.Assumption.on.A}
 I_m(A)=J_1\cdot J_2=J_1\cap J_2\sset R, \text{ each $J_i$ contains a regular element},\quad\quad       I_{m-1}(A)\sseteq J_1+J_2.
\eeq
  \eee
\noindent We prove: if these conditions hold for $A\in \Mat{R}$ then they also hold for $\bA\in \Mat{S}$.

This reduces the initial decomposability question (over $R$) to that over $S$. The study goes via the ``obstruction to decomposability",
 as explained in \S\ref{Sec.Intro.Structure.Contents}.

\subsection{Reduction to Noetherian rings}\label{Sec.Reduction.to.Noetherian}
\bprop
Let $A\in \Mat{R}$.  Then   $I_m(A)=J_1\cdots J_r$ and $A$ is $(J_1,\dots, J_r)$-decomposable \iff there exists a Noetherian subring $S\sseteq R$
 that contains the entries of $A$, with $I_m(A)=\prod \bJ_i\sset S$, such that $A$ is $(\bJ_1,\dots, \bJ_r)$-decomposable over $S$.
\hfill (The ideals $\bJ_i\sset S$ are defined in the proof.)
\eprop

Similar statement holds for the stable-decomposability.
\bpr The part $\Lleftarrow$ is trivial. We prove the part $\Rrightarrow$.

  Assume $U\cdot A\cdot V^{-1}=\oplus A_i$ for some $U\in GL(m,R)$, $V\in GL(n,R)$ and $A_i\in \operatorname{Mat}_{m_i\times n_i}(R)$.
 Take the $\Z$-subalgebra $S\sseteq R$ finitely generated by $1\in R$ and by the entries of $A$, $\{A_i\}$, $U$, $U^{-1}$, $V$, $V^{-1}$.
 We keep all the polynomial relations (with coefficients in $S$) among these elements that hold in $R$. Thus $S$ is finitely generated, in particular Noetherian.

The ideal  $J_i\!=\!I_{m_i}(A_i)\!\sset\! R$ is (finitely)
  generated by the maximal minors of $A_i$. Define $\bJ_i\!\sset\! S$ by these same maximal minors. Then
   $U\!\cdot\! A\!\cdot\! V^{-1}\!=\!\oplus A_i$, for $U\!\in\! GL(m,S)$, $V\!\in\! GL(n,S)$, and $I_{m_i}(A)\!=\!\bJ_i.$
\epr

\

 For square matrices, $A\in \Mats{R}$, the conditions \eqref{Eq.Assumption.on.A}  read:
\beq\label{Eq.Assumption.on.A.square}
 \det(A)=f_1\cdot f_2 \text{ is regular in }R,\quad\quad (f_1)\cdot (f_2)=(f_1)\cap (f_2)\sset R,\quad\quad     I_{n-1}(A)\sseteq (f_1,f_2).
\eeq

 We can often pass to a Noetherian subring $S\sseteq R$ while preserving these conditions:
\bprop\label{Thm.reduction.to.Noetherian}
Assume the conditions \eqref{Eq.Assumption.on.A.square} hold  for $A\in \Mats{R}$. Assume for any finitely generated subring $S\sseteq R$, such that
 $A\in \Mats{S}$, the chain of ring extensions
\beq
S\sseteq S:_R (f_1f_2)\sseteq S:_R (f_1f_2)^2\sseteq\cdots=\liml_\to \big(S:_R(f_1f_2)^\bullet\big)
\eeq
 stabilizes.
 Then there exists a Noetherian subring $S\sseteq R$ such that $A\in \Mats{S}$, and  the conditions \eqref{Eq.Assumption.on.A.square} hold  over $S$.
\eprop
\bpr
\bee[\bf Step 1.\!]
\item Take the entries $\{a_{ij}\}$ of $A$.
 For each $(n-1)$-block  expand  $\det(A_\Box)=  d^{(1)}_{\Box}\cdot f_1+  d^{(2)}_{\Box }\cdot  f_2$.
  Take the subring $S\sseteq R$ generated by the (finite) collection $\{a_{ij}\}$, $f_1,f_2$,
   $\{d^{(1)}_{\Box}\}$, $\{d^{(2)}_{\Box}\}$.
     We have:
 \beq\label{Eq.Inner.Conditions}
A\in \Mats{S},\quad\quad \det(A)=f_1\cdot f_2,\quad  I_{n-1}(A)\sseteq S\cdot f_1+S\cdot f_2\sseteq S.
 \eeq
The ring $S$ is Noetherian, being a finitely-generated $\Z$-algebra.

 But the condition  $(S\cdot f_1)\cdot (S\cdot f_2)=(S\cdot f_1)\cap (S\cdot f_2)$ does not necessarily hold,
 as $S$ is not necessarily a UFD, and  we do not have the submodule/ideal contraction property,
 and cannot use $S\cdot f_i=(Rf_i)\cap S$.

\item
Note that  conditions \eqref{Eq.Inner.Conditions} hold also over any further extensions of $S$, i.e.,  if $S\sset\tS$ then
$\det(\tS\otimes A)=\tS\cdot (f_1f_2)\sset \tS$ and $I_{n-1}(\tS\otimes A)\sseteq \tS f_1+\tS f_2.$

 We claim:  $S f_1\cap S f_2=(f_1f_2)\cdot (S:_R(f_1,f_2))$. Indeed, any element of $S f_1\cap S f_2$ is presentable as $f_1f_2 c$,
 for some $ c\in R$ satisfying $f_1 c,f_2 c\in S$. And vice-versa.

Thus we want an extension $S\sseteq S'\sseteq R$ satisfying: $S'$ is Noetherian and $S':_R (f_1,f_2)=S'$.
 Define $S':=\sum^\infty_{j=0} S:_R(f_1\cdot f_2)^j$ and observe:
 \beq
 S'\sseteq S':_R (f_1,f_2) \sseteq S':_R (f_1f_2)=S'.
 \eeq
Moreover, by our assumption the sum $\sum_d S:_R(f_1\cdot f_2)^d$ stabilizes.
 Therefore, to show that $S'$ is finitely generated (hence Noetherian), it is enough to verify: the subring $ S:_R(f_1\cdot f_2)\sset R$ is finitely generated.
 Indeed,  the ideal $(f_1f_2)\cdot (S:_R(f_1 f_2))\sset S$  is finitely generated, as $S$ is Noetherian. And $f_1\cdot f_2\in R$ is a non-zero divisor.
  Thus $S:_R(f_1,f_2)\sset R$ is finitely generated.
\epr
\eee
\beR
\bee[\bf i.]
\item The assumption ``the chain  $S\sseteq S:_R (f_1f_2)\sseteq S:_R (f_1f_2)^2\sseteq\cdots$ stabilizes"
 is satisfied if the Artin-Rees condition (for any regular element $f\in R$) holds: $S\cap (f^N\cdot R)\sset f^{N-d}S$ for $N\gg1$.
 And this holds for numerous non-Noetherian rings.

\item This assumption  is non-empty. For example,
 let $R=\k[t,\frac{x_1}{t},\frac{x_2}{t^2},\dots]$ and $S=\k[t]$. Then $\sum S:_R t^d=R$, non-Noetherian, and the chain $\{S:_R t^\bullet\}$ does not stabilize.

Another example is $R=\k[x,y]\supset S=\k[xy]$. Here $\sum S:_R y^d=R$, but the chain $\{S:_R y^\bullet\}$ does not stabilize.
\eee
\eeR

\subsection{The obstruction to decomposability}\label{Sec.Obstruction.To.Decomposability}
 Let  $A\in \Mat{R}$ with  $I_m(A)=J_1\cdot J_2=J_1\cap J_2\sset R.$
 We construct the obstruction to decomposability, as is explained in \S\ref{Sec.Intro.Structure.Contents}. Let $M:=\Coker(A)$.

\bprop\label{Thm.Obstruction.to.Decomposability}
The exact sequence  of    $R$-modules  (constructed in the proof)
\beq\label{Eq.Exact.Sequence.For.Decomposability}
0\to H^0_{J_1+J_2}(M)\to  M\to M_1\oplus M_2\to Q\to0
\eeq           satisfies:
        \bee[\bf 1.]
\item The modules $M_1,\!M_2,\!Q$ are    finitely generated and satisfy: $J_i\!\cdot\! M_i\!=\!0$, $(J_1\!+\!J_2)\!\cdot\! Q\!=\!0,$ $Ann(M_i)\!\sseteq\! \sqrt{J_i}.$
\item If $M$ is $(J_1,J_2)$-decomposable then $Q=0$.
\item  If $Q\!=\!0$ and the submodule $H^0_{J_1+J_2}(M)\sset M$ splits off as a direct summand, then $M$ is $(J_1,J_2)$-decomposable.
\eee
     \eprop
     \bpr      Separate the components of the support: \quad $\amalg \Spec(\quots{R}{J_i})  \stackrel{\pi}{\to }\Spec(\quots{R}{ \prod J_i })$.
  Here $\pi$ is defined via the idempotents,  $e_i\in \quots{R}{J_i}$:
  \beq\label{Eq.Separating.Components}
  \quots{R}{I_m(A)}\stackrel{\pi^*}{\to} \   \quots{R}{J_1}\times\quots{R}{J_2},\quad\quad x\to (e_1\cdot x,e_2\cdot x).
  \eeq
Note that $\pi^*$ is injective. Indeed, $Ker(\pi^*)=   J_1\cap J_2=J_1\cdot J_2=0\sset \quots{R}{ J_1 J_2 }$.

Restrict the module $M$ to $V(J_i)$, i.e.,  define $M'_i:=\Coker(A\otimes \quots{R}{J_i})$. Kill the unwanted torsion,
\beq\label{Eq.Lifted.M.killed.torsion}
0\to H^0_{J_1+J_2}(M'_i)=H^0_{J_j}(M'_i)\to M'_i\to M_i\to 0, \quad \text{ for $i\neq j$}.
\eeq
 Now push the $\quots{R}{J_i}$-modules $M_i$  back to $Spec(R)$, i.e.,  consider these as $R$-modules.

\bee[\bf\text{Part }1.]
\item Equation \eqref{Eq.Separating.Components} gives the morphism of $R$-modules $M\stackrel{\oplus\phi_i}{\to} \oplus M_i$.
 Here $ker(\phi_i)=H^0_{J_j}(M)$ for $j\neq i$. Therefore $ker(\phi_1\oplus\phi_2)=H^0_{J_1+J_2}(M)$.
 Hence the sequence \eqref{Eq.Exact.Sequence.For.Decomposability} is exact at $M$. By construction we have $J_i\cdot M_i=0$.

We prove: $Ann(M_i)\sseteq\sqrt{J_i}\sset R.$ Indeed, $Ann(M_i)\cdot M'_i\sseteq H^0_{J_1+J_2}(M'_i)$ by \eqref{Eq.Lifted.M.killed.torsion}.
  Thus $(J_1+J_2)^d\cdot Ann(M_i)\cdot M'_i=0$ for some $d\gg1.$ In particular, $J_j^d\!\cdot\! Ann(M_i)\!\cdot\! M'_i=0,$ i.e.,
    $J_j^d\!\cdot \! Ann(M_i)\!\cdot \!M\!\sseteq \!J_i\!\cdot \!M.$ But then $J_j^{d+1}\cdot Ann(M_i)\cdot M=0.$
     Thus $J_j^{d+1}\!\cdot\! Ann(M_i)\!\sseteq\! Ann(M)\!\sseteq\!\sqrt{J_1\!\cdot\! J_2}.$ And then, by checking the minimal associated primes, we get  $Ann(M_i)\!\sseteq\!\sqrt{J_i}.$

\medskip

The module $Q$ is defined as the quotient by  \eqref{Eq.Exact.Sequence.For.Decomposability}. Thus $M_1,M_2,Q$ are finitely generated.

\medskip

 We verify: $(J_1+J_2)\cdot Q=0$.
 Fix an element of $Q$ and take its representative $z_1\oplus z_2\in M_1\oplus M_2$. As  $\phi_1,\phi_2$ are surjective, we have: $z_i=\phi_i(\tz_i)$
  for some $\tz_1,\tz_2\in M$. Therefore:
  \beq
  M_1\oplus M_2\supset (J_1+J_2)(z_1\oplus z_2)=J_2z_1\oplus J_1z_2=J_2\cdot \phi_1(\tz_1) \oplus J_1\cdot \phi_2(\tz_2)=
  \eeq
  \[
  =   J_2(\phi_1(\tz_1)\oplus\phi_2(\tz_2))+J_1(\phi_1(\tz_1)\oplus\phi_2(\tz_2))\in (J_1+J_2)\cdot Im(\phi_1\oplus\phi_2).\]
   Therefore $ (J_1+J_2)(z_1\oplus z_2)$ goes to   $0\in Q$,  by the exactness of \eqref{Eq.Exact.Sequence.For.Decomposability}.

\item  Suppose $M=N_1\oplus N_2$ with $\sqrt{Fitt_0(N_i)}=\sqrt{J_i}$. Equation \eqref{Eq.Lifted.M.killed.torsion} reads then:
\beq
0\to (\quots{R}{J_i}\otimes  N_j)\oplus  H^0_{J_j}(\quots{R}{J_i}\otimes  N_i)\stackrel{\psi}{\to} \quots{R}{J_i}\otimes  (N_1\oplus N_2)\to  M_i\to0
 \quad \text{ for }i\neq j.
\eeq
 Here $\psi$ is block-diagonal, with one of the blocks   $\quots{R}{J_i}\otimes  N_j\isom{Id}  \quots{R}{J_i}\otimes  N_j .$
  This gives: $M_i=\quots{\quots{R}{J_i}\otimes  N_j}{H^0_{J_j}(\quots{R}{J_i}\otimes  N_i)}$.
 In particular, $M\twoheadrightarrow M_1\oplus M_2.$ Hence $Q=0.$

\item
 If $Q\!=\!0$ (and the submodule $H^0_{J_1+J_2}\!(M)\!\sset\! M$ splits off) then $M\!=\!H^0_{J_1+J_2}\!(M)\!\oplus\! \!M_1\!\oplus \!M_2.$
  \\Here $\sqrt{\!Ann(H^0_{J_1\!+\!J_2})}\!\supseteq\! J_1\!+\!J_2.$ Therefore \!$\sqrt{\!Fitt_o[Ann(H^0_{J_1+J_2})\!\oplus\! M_1]}\!=\!\sqrt{\!J_1}$ and
  \!$\sqrt{\!Fitt_o( M_2)}\!=\!\sqrt{\!J_2}.$
   Hence the claimed decomposition.
\epr
\eee

\subsection{Functoriality of the obstruction $\boldsymbol Q$} Given a morphism of rings $R\to S$, take the images
\beq
\bA\in \Mat{S}, \qquad  I_m(\bA)\!=\!\bJ_1\!\cdot\! \bJ_2\!\sset\! S, \qquad
 \bbM\!:=\!\Coker(\bA).
 \eeq
 As in \S\ref{Sec.Obstruction.To.Decomposability} we restrict $M$ to $V(J_i)\!\sset\! \Spec(R)$  and take the torsion on $V(J_1\!+\!J_2)$, i.e.
  $H^0_{J_1\!+\!J_2}(\quots{R}{J_i}\! \otimes\! M)\!\sset\! \quots{R}{J_i} \! \otimes\! M$.
 The analogue of the sequence \eqref{Eq.Exact.Sequence.For.Decomposability} over $S$ can be obtained in two ways: either by applying $S\otimes,$ or starting directly from $\bbM.$ Accordingly we get the modules $S\otimes Q$ and $\bQ.$
  Furthermore, we can restrict $\bbM$
 to $V(\bJ_i)\!\sset\! \Spec(S)$ in two ways: either by
  $S\!\otimes\! H^0_{J_1+J_2}(\quots{R}{J_i} \!\otimes M)\!\to \!S\!\cdot \! H^0_{J_1+J_2}(\quots{R}{J_i} \!\otimes \!M)\!\sseteq \!\quots{S}{\bJ_i} \!\otimes\!\bbM$
  or as $ H^0_{\bJ_1\!+\!\bJ_2}(\quots{S}{\bJ_i} \!\otimes\!\bbM)\!\sseteq \!\quots{S}{\bJ_i} \!\otimes\!\bbM.$

\bprop\label{Thm.Obstruction.Functoriality}
 Suppose $\bJ_1\cdot \bJ_2=\bJ_1\cap \bJ_2\sset S$ and let $i\in 1,2$. Then:
\bee[\bf 1.]
\item The sequence \eqref{Eq.Exact.Sequence.For.Decomposability} induces the commutative diagram:
\beq\label{Eq.Obstruction.Base.Change.Diagram}
\bM
S\otimes M&\to & \oplus_i (S\otimes M_i)&\to& S\otimes Q&\to&0
\\
\rotatebox[]{-90}{$\isom{}$}&&\de_0\rotatebox[]{-90}{$\twoheadrightarrow$}&&\phi\rotatebox[]{-90}{$\twoheadrightarrow$}
\\
\bbM&\to&\oplus_i \bbM_i&\to &\bQ&\to&0
\eM
\eeq
Here $\de_o$ and $\phi$ are surjective.
\item The (natural) morphism $S\cdot H^0_{J_1+J_2}(\quots{R}{J_i}\otimes M)\to H^0_{\bJ_1+\bJ_2}(\quots{S}{\bJ_i}\otimes \bbM)$ is an embedding.
 It is an isomorphism (for both $i=1$ and $i=2$) \iff $\phi:\ S\otimes Q\to \bQ$ is an isomorphism.
\eee
\eprop
\bpr
\bee[\bf 1.]
\item One can either apply $S\otimes$ to the sequence \eqref{Eq.Lifted.M.killed.torsion} or   write the corresponding sequence for the module $\bbM$.
 This gives the two rows of the diagram:
\beq\label{Eq.Obstruction.Functoriality}
\bM
&S\cdot  H^0_{J_1+J_2}(M'_i)&\stackrel{\ep}{\hookrightarrow} &   S\otimes M'_i&\to& S\otimes M_i&\to&0
\\
&\de_2\rotatebox[]{-90}{$\hookrightarrow$}&&\de_1\rotatebox[]{-90}{$\isom{}$}&&\de_0\rotatebox[]{-90}{$\twoheadrightarrow$}
\\
0\to& H^0_{\bJ_1+\bJ_2}(\bbM'_i)&\to& \bbM'_i&\to &\bbM_i&\to&0
\eM
\eeq
Here the isomorphism $\de_1$ is the natural composition:
\beq
S\otimes M'_i=S\otimes\Coker(A\otimes \quots{R}{J_i})\isom{}\Coker(A\otimes \quots{R}{J_i}\otimes S)=\Coker(\bA \otimes \quots{S}{\bJ_i})=:\bbM'_i.
\eeq
To define $\de_2$   we observe that the morphisms $\de_1,\ep$ in \eqref{Eq.Obstruction.Functoriality} are $S$-linear. Therefore for any $\xi\in S\cdot  H^0_{J_1+J_2}(M'_i)$ and a corresponding $d\gg1$ one has:
 $(\bJ_1+\bJ_2)^d\cdot \de_1(\ep(\xi))=0$. Thus $(\de_1\circ\ep)(S\cdot  H^0_{J_1+J_2}(M'_i))\sseteq H^0_{\bJ_1+\bJ_2}(\bbM'_i)$.
  Hence the morphism $\de_1\circ\ep$ factorizes through $H^0_{\bJ_1+\bJ_2}(\bbM'_i)$, thus defining $\de_2$.
  Moreover, $\de_2$ is injective, as $\de_1,\ep$ are injective.

  Finally, $\de_0$ is defined (and is surjective) by  the standard diagram chasing.

\

Now the surjectivity of $\de_0$ gives the claimed diagram \eqref{Eq.Obstruction.Base.Change.Diagram}. The existence/surjectivity of $\phi$ follows by
 the standard diagram chasing.
\item
By the diagram chasing in \eqref{Eq.Obstruction.Functoriality} we get: $\de_2$ is an isomorphism \iff $\de_0$ is an   isomorphism \iff $\phi$ is an isomorphism.
\epr
\eee

\subsection{Reduction of decomposability question to local/Henselian/complete  rings}\label{Sec.Reduction.to.local.complete}
\bcor\label{Thm.reduction.to.complete.local.ring.1}
Let $A\in \Mat{R}$ with  $I_m(A)\!=\! J_1\cdot J_2\!=J_1\cap J_2$. Assume  the submodule $H^0_{J_1+J_2}(M)\sset M$ splits off as a direct summand.
\bee[\bf 1.]
\item $A$ is stably-$(J_1, J_2)$-decomposable \iff for   every localization
 $R \to  R_\cm$ at a  maximal ideal $J_1+J_2\sseteq \cm\sset R$
  the matrix   $\bA\in \Mat{R_\cm}$ is $(\bJ_1   ,\bJ_2  )$-decomposable.
\item Let $(R,\cm)$  be a local ring and take the   Henselization   $R \to  R^{h}$.
  Then $A$ is  $(J_1, J_2)$-decomposable \iff
  $\bA \in \Mat{R^{h}}$ is $(\bJ_1 ,\bJ_2 )$-decomposable.
\item Let $(R,\cm)$  be a local ring and assume that the completion functor $\otimes \hR^{(\cm)}$
 is exact on finitely generated $R$-modules.
 Then $A$ is $(J_1, J_2)$-decomposable \iff
  $\bA\in \Mat{\hR^{(\cm)}}$ is $(\bJ_1   ,\bJ_2 )$-decomposable.
\eee\ecor
\bpr As  the submodule $H^0_{J_1+J_2}(M)\sset M$ splits off, the only obstruction to decomposability is $Q.$
   The functors of
 localization,  $\otimes R_\cm$, Henselization, $\otimes R^{(h)}$, and (by our assumption) completion, $\otimes \hR^{(\cm)}$,
  are exact.  Hence they preserve the exactness of the sequence \eqref{Eq.Exact.Sequence.For.Decomposability}.
  \\Then part 2 of proposition \ref{Thm.Obstruction.Functoriality} gives: $S\!\otimes\! Q\!\cong\! \bQ$. Finally:
\bee
\item The module $Q$ vanishes \iff   its  localizations at all the maximal ideals vanish. And it is enough to localize at the ideals satisfying $J_1+J_2\sseteq \cm$, as
   $Q_\cm=0$ for $\cm\not\supseteq J_1+J_2$.
\item  The module $Q$ over a local ring vanishes \iff   its $\cm$-Henselization vanishes.
\item
  The (finitely generated) module $Q$ over a local ring $(R,\cm)$ vanishes \iff $  \hR^{(\cm)}\otimes Q=0$.
\eee
\epr

\beR
  The completion functor $\otimes\hR^{(\cm)}$ is exact on finitely generated modules over local Noetherian rings. But for
 local non-Noetherian rings  $\otimes\hR^{(\cm)}$ is not always exact, see  example 2.8.7 in \cite{Schenzel-Simon}.
\eeR

 We restate the conclusion for modules.
 \bcor Take a finitely presented $R$-module  $M$ with factorized Fitting ideal $Fitt_0(M)=J_1\cdot J_2=J_1\cap J_2$.
  Suppose  the submodule $H^0_{J_1+J_2}(M)\sset M$ splits as a direct summand. TFAE:
    \bee[\bf i.]
    \item $M$ is $(J_1,J_2)$-decomposable;
    \item  the localization $M_\cm$ is $((J_1)_\cm,(J_2)_\cm)$-decomposable
  for any maximal ideal   $\cm\supseteq J_1+J_2$
  \item    (for $(R,\cm)$, assuming exactness of  $-\otimes\hR^{(\cm)}$)  the completion $M\otimes\hR^{(\cm)}$ is decomposable.
   \eee
\ecor

\bex (Semi-local rings/matrices over multi-germs of spaces)
 Let $R=R_1\times\cdots\times R_k$, where $\{(R_l,\cm_l)\}$ are local rings. Fix the corresponding idempotents,
\beq
1=\sum^k_{l=1} e_l,\quad\quad e_l\cdot e_\tl=\de_{l,\tl}\cdot e_l.
\eeq
 Let $A\in\Mat{R}$ with $I_m(A)= J_1J_2=J_1\cap J_2.$ Suppose $H^0_{J_1+J_2}(Coker(A))=0.$
 Corollary \ref{Thm.reduction.to.complete.local.ring.1} gives:
   $A$ is stably $(J_1,J_2)$-decomposable \iff  $e_l A\in \Mat{R_l}$ is $(e_l J_1,  e_l J_2)$-decomposable,
   $e_l A\sim   A^{(l)}_1\oplus A^{(l)}_2$, for each $l=1\dots k$. We do not assume that all $\{A^{(l)}_1\}_l$  (or all $\{A^{(l)}_2\}_l$) are of the same size.
    In fact, by \eqref{Eq.Decompos.of.A.vs.Coker(A)}  (and also directly via the idempotents) one gets:
$A\stackrel{s}{\sim} \oplus^k_{l=1}  e_l A$.
Assuming the decomposability for each $l$, i.e.,  $ e_l A\sim \oplus^2_{i=1}A^{(l)}_i$, we get:
$A\stackrel{s}{\sim} \oplus^2_{i=1}  \oplus^k_{l=1}   A^{(l)}_i$.

Geometrically, for $\Spec(R)=\amalg \Spec(R_l)$ (the finite union of germs of spaces), $A$ is stably decomposable \iff each restriction $\{A|_{\Spec(R_l)}\}$
 is decomposable.
\eex

\subsection{The behaviour of the conditions \eqref{Eq.Assumption.on.A} under   base change}
\bel
Assume  the  conditions \eqref{Eq.Assumption.on.A} hold for $A\in \Mat{R}$. Then they hold for localizations, $\bA\in \Mat{R_\cm}$, Henselizations,
 $\bA\in \Mat{R^h}$, and (when the functor $\otimes\hR^\cm$ is exact) completions, $\bA\in \Mat{\hR^{\cm}}$.
\eel
\bpr The persistence of factorization, $I_m(A)=J_1\cdot J_2$, and inclusion, $I_{m-1}(A)\sseteq J_1+J_2$, is obvious. If $x\in J_i$ is regular
 then $\bx\in\bJ_i$ is regular. (Apply the needed flat base change to the exact sequence $0\to R\stackrel{\times x}{\to}R\to \quots{R}{(x)}\to0$.)
  For the condition $J_1\cap J_2=J_1\cdot J_2$ use its equivalent form, lemma \ref{Thm.coprime.ideals}, and observe that $Tor_\bullet$
   is functorial under flat base-change.
\epr

\subsection{Reduction to completion for $C^\infty$-rings}\label{Sec.Reduction.Cinfty.to.complete}
 Let $R$ be one of the standard (non-Noetherian) rings of Analysis:
\bei
\item
smooth functions on an open set, $C^\infty(\cU)$, for $\cU\sseteq \R^p$;
\item
germs of smooth functions along a closed subset, $C^\infty(\cU,Z)$;
\item
 $\quots{C^\infty(\cU)}{I}$,  $\quots{C^\infty(\cU,Z)}{I}$,  for an ideal $I$, i.e.,  the rings of functions on the ``subscheme" $V(I)\sset  \cU $.
\eei
In this case the completion functor $\otimes\hR$ is far  from being exact/faithful.
   For example, the sequence $[0\to (f)\to R]\otimes\hR$ is not exact for any $0\neq f\in \cm^\infty,$ and $\cm^\infty\otimes\hR=0.$
 While our main decomposability criterion (theorem \ref{Thm.Decomposition.Square.Matrices}) is
 applicable to such rings (see \S\ref{Sec.Decomposition.Remarks}.iv), the ``reduction-to-completion" is still important.
 E.g., smooth functions are often studied via their Taylor expansions. Below we give an independent reduction-to-completion criterion.

Take a closed subset $Z\!\sset\! \cU$,
 and the corresponding completion $\hR^{(Z)}\!:=\!\liml_{\leftarrow}\quots{R}{I(Z)^{\bl j\br}}$. Here the differential power of ideal, $I(Z)^{\bl j\br},$
consists of functions vanishing at all the points of $Z$ to the order at least $j$:
 \beq
 I(Z)^{\bl j\br}:=\cap_{x\in Z}\cm^j_x\sset R.
\eeq
The kernel of the map $R\!\to\! \hR^{(Z)}$ is the ideal  $I(Z)^{\bl \infty\br}\!:=\!\cap_{x\in Z} \cm^\infty_x$ of functions that are flat on $Z.$

The completion at $Z$ erases the complement of the formal neighbourhood of $Z$, therefore we always assume $V(I_m(A))\sseteq  Z\sset\cU$.
 Just a slight strengthening of this condition binds the decomposability over $R$ to that over $\hR^{(Z)}$.
  Let the ring $R$ be  $\quots{C^\infty(\cU)}{I}$ or  $\quots{C^\infty(\cU,Z)}{I}$.
Take the completion $R\to \hR^{(Z)}$ and accordingly $\Mat{R}\stackrel{\phi}{\to} \Mat{\hR^{(Z)}}$.
\bel
 Assume $I_m(A)\supseteq I(Z)^{\bl \infty\br}$.
 Then $A\sim\oplus A_i$ \iff $\phi(A)\sim \oplus \hA_i$. (And then $\hA_i\sim \phi(A_i)$.)
\eel
\bpr (the direction $\!\Lleftarrow$)
Whitney's extension theorem ensures the surjectivity of completion, $R\stackrel{\phi}{\twoheadrightarrow}\hR^{(Z)}$,
 see   e.g., \S1.5 of \cite{Narasimhan} or \cite[\S2]{Belitski.Kerner.Whitney}.
 Therefore we choose $R$-representatives $\{A_i\}$ of $\{\hA_i\}.$ W.l.o.g. we assume: the matrices $\{A_i\}_i$ are of full ranks on $\cU\smin Z.$
 We have:
 $\phi(A)=\hU_m\cdot (\oplus \phi(A_i))\cdot \hU^{-1}_n$, for some $\hU_m\in GL(m,\hR^{(Z)})$,
  $\hU_n\in GL(n,\hR^{(Z)})$.

We claim the surjectivity of group-homomorphism  $GL(n,R)\twoheadrightarrow GL(n,\hR^{(Z)})$. Indeed, take some Whitney representative $U_n\in \Mats{R}$ of $\hU_n\in GL(n,\hR^{(Z)})$,
 then $V(det(U_n))\cap Z=\empty$. Thus $U_n$ is invertible in a small neighbourhood of $Z$. Modify $U_n$ outside of this neighbourhood (e.g.
  by   cutoff functions) to achieve the invertibility on the whole $\cU$.

Finally, take some representatives $U_m\in GL(m,R)$, $U_n\in GL(n,R)$ of $\hU_m,\hU_n$. We get:
\beq
A= U_m\cdot(\oplus A_i)\cdot U^{-1}_n+A^\infty,\quad
 \text{ where }A^\infty\in \Mat{I(Z)^{\bl \infty\br}}.
 \eeq
 Moreover, we can ssume: $A^\infty=\zero$ far from $Z.$ (If needed, one adjusts $U_m,U_n$.)

  As $I_m(A)\supseteq I(Z)^{\bl \infty\br}$, we can present $A^\infty=A\cdot C$, for some $C\in \Mats{R}$.
  Moreover, we can assume $C\in \Mats{I(Z)}$. Indeed, take the $l_2$-matrix norm $\Mat{\R}\stackrel{\|-\|}{\to}\R.$
  Then
    $A^\infty\cdot \frac{1}{\sqrt{||A^\infty||}}\in \Mat{I(Z)^{\bl \infty\br}}$.  Thus $A^\infty\cdot \frac{1}{\sqrt{||A^\infty||}}=A\cdot \tC$,
     hence  $A^\infty =A\cdot \tC\cdot \sqrt{||A^\infty||}.$ Finally, $C=\zero$ far from $Z.$

Therefore the matrix $\one\!-\!C$ is invertible on $\cU$.
    We get:  $A\!=\! U_m\!\cdot\!(\oplus A_i)\!\cdot\! U^{-1}_n\!\cdot\!(\one-C)^{-1}.$
\epr

\bex
Let $R=C^\infty(\R^p)$ and assume $I_m(A)\supseteq (\ux)^\infty$, the ideal of functions flat at $o\in \R^p$. Then $A$ is decomposable \iff its Taylor expansion at
 $o$  is decomposable.
\eex

We recall the standard way to ensure the assumption  $I_m(A)\supseteq I(Z)^{\bl \infty\br}$ of the last lemma.
\bel
Let $R=C^\infty(\cU)$,  $A\in \Mat{R}$. Assume the function $det(AA^T): \cU \to \R$ satisfies the
  {\L}ojasiewicz-type inequality: $|\det(AA^T)(x)|\ge C\cdot dist(x,Z)^\de$, for some constants $C,\de>0$ and any $x\in \cU$.
 Then  $I_m(A)\supseteq I(Z)^{\bl \infty\br}$.
\eel
\bpr
 Recall that $I_m(A)\ni \det(A  A^T),$ by Cauchy-Binet formula. Therefore for any $h\in I(Z)^{\bl \infty\br}$ it is enough to verify: $\frac{h}{det(AA^T)}\in R$.
 By our assumptions, all the  derivatives of  $\frac{h}{det(AA^T)}\in R$  tend to 0 on $Z$. Therefore this ratio extends to a flat $C^\infty$ function on $Z$.
\epr

\subsection{Passage to the quotient ring, the case of non-flat base change}\label{Sec.Reduction.Passage.to.Quotient}

In many cases the (in)decomposability is preserved under the base-change even when the morphism $R\to S$ is non-flat.
\bprop\label{Thm.Decomposability.Restriction.to.sublocus} Suppose the functor $S\otimes$ is faithful on finitely generated $R$-modules.
Suppose:
\bee[\bf i.]
\item  $I_m(A)=J_1\cdot J_2=J_1\cap J_2\sset R$ and $\bJ_1\cdot \bJ_2=\bJ_1\cap \bJ_2\sset S$, see \S\ref{Sec.Background.Notations.Conventions}.ii.
\item the submodule $H^0_{J_1+J_2}(M)\sset M$ splits off,  
 and moreover
  $S\cdot H^0_{J_j}( M)= H^0_{\bJ_j}( \bbM)$ for $j=1,2$.
  \eee
 Then $M$ is $(J_1,J_2)$-decomposable \iff $\bbM$ is $(\bJ_1,\bJ_2)$-decomposable.
\eprop
\bpr (of the direction $\Lleftarrow$)
We have $\bQ=0$ and want to deduce $Q=0.$ (Then one invokes proposition \ref{Thm.Obstruction.to.Decomposability}.)
 By part 2 of proposition \ref{Thm.Obstruction.Functoriality} it is enough to verify:
 the embedding $S\cdot H^0_{J_1+J_2}(\quots{R}{J_i}\otimes M )\hookrightarrow  H^0_{\bJ_1+\bJ_2}(\quots{S}{\bJ_i}\otimes \bbM )$ is an isomorphism.

For $i\neq j$ take the submodule $\quots{R}{J_i}\cdot H^0_{J_j}( M)\sseteq \quots{R}{J_i}\otimes M  ,$ see \S\ref{Sec.Background.Notations.Conventions}.ii.
 Obviously  $\quots{R}{J_i}\cdot H^0_{J_j}( M)\sseteq   H^0_{J_j}( \quots{R}{J_i}\otimes M).$
 We claim: this is an equality.
    Indeed, take an element of $H^0_{J_j}(\quots{R}{J_i}\otimes M)$ and let $z\in M$ be its representative. Then
    $J^d_j\cdot z\sseteq J_iM\cap J_jM$ for some $d\gg1$. But then  $J^{d+1}_j\cdot z=0\sset M.$
     Thus $z\in H^0_{J_j}( M)$, proving the claim.

    In the same way one gets:
 $\quots{S}{\bJ_i}\cdot H^0_{\bJ_j}( \bbM)= H^0_{\bJ_j}(\quots{S}{\bJ_i}\otimes \bbM)\sseteq \quots{S}{\bJ_i}\otimes \bbM.$
Altogether:
\beq
S\cdot H^0_{J_1+J_2}(\quots{R}{J_i}\otimes M)=S\cdot H^0_{J_j}(\quots{R}{J_i}\otimes M)=
(S\cdot \quots{R}{J_i})\cdot H^0_{J_j}( M)\stackrel{!}{=}\quots{S}{\bJ_i}\cdot H^0_{\bJ_j}( \bbM)=
\eeq
\[\hspace{8cm}
= H^0_{\bJ_j}(\quots{S}{\bJ_i}\otimes \bbM)=
  H^0_{\bJ_1+\bJ_2}(\quots{S}{\bJ_i}\otimes \bbM).
\]
(The equality  `!' holds by the assumption ii.)

Part 2 of proposition \ref{Thm.Obstruction.Functoriality} gives: the morphism
 $\phi:\ S\otimes Q\to \bQ$ is an isomorphism. By our assumption $\bQ=0$, thus $S\otimes Q=0$. By our assumption  the functor $S\otimes $ is faithful,
  therefore $Q=0$. Finally (proposition \ref{Thm.Obstruction.to.Decomposability}): $M$ is $(J_1,J_2)$-decomposable.
\epr

\bcor\label{Thm.Decompos.over.R.vs.R/a.Square.Case}
(For square matrices) Let $S=\quots{R}{\ca}$ and  $A\in \Mats{R}$ with $det(A)=f_1\cdot f_2$. Suppose:
\bee[\bf i.]
\item $f_1,f_2\in R$ are regular, coprime, and $ \bar{f}_1,  \bar{f}_2\in S$ are regular, coprime.
 \item  $S\cdot (Im(A):f_j)=Im(\bA):\bar{f}_j\sset S^n$ for $j=1,2$.
 \eee
   Then $A$ is $f_1,f_2$-decomposable \iff $\bA$ is $\bar{f}_1,\bar{f}_2$-decomposable.
\ecor
\bpr By lemma \ref{Thm.coprime.ideals}: $H^0_{f_1,f_2}(M)=0$ and  $H^0_{\bar{f}_1,\bar{f}_2}(\bbM)=0.$
 The assumption ii. implies   $S\cdot H^0_{(f_j)}( M)= H^0_{(\bar{f_j})}( \bbM)$.
 Now apply proposition \ref{Thm.Decomposability.Restriction.to.sublocus}.
\epr

\beR
The assumption  $S\cdot (Im(A):f_j)=Im(\bA):\bar{f}_j$ is needed here. For example, let $A=\bbm x&y\\0&z\ebm$, and $S=\quots{R}{(y)}$.
 Then $\bA$ is $\bx,\bz$-decomposable, though $A$ is indecomposable.
\eeR

\subsection{Graded-to-local reduction}\label{Sec.Reduction.Graded-to-Local}
 Take an $\N$-graded ring, $R=\oplus_{d\in \N} R_d$, with $R_0$ local, Noetherian. Denote by $\cm\sset R$ the ideal generated by $R_{>0}$
   and by the maximal ideal of $R_0$. Thus $\cm$ is a homogeneous maximal ideal. It is also the largest among the homogeneous ideals.
    The localization $R\to R_\cm$  is injective.
\bel\label{Thm.Decomposability.Graded-to-Local.reduction}
Assume $A\in \Mat{R}$  is graded (\S\ref{Sec.Background.Notations.Conventions}.v).
 Then $A\sim \oplus A_i$ by  $GL(m,R)\times GL(n,R)$ \iff 
  $A\sim \oplus A_i$ by  $GL(m,R_\cm)\times GL(n,R_\cm)$.
\eel
\bpr (of the part $\Lleftarrow$)
  One has $\Coker(A)\otimes R_\cm\cong \Coker(\oplus A_i)\otimes R_\cm$. Recall: two graded modules are
 isomorphic (in the non-graded way, i.e. up to the shifts of grades) \iff their localizations at the maximal among the homogeneous ideals are isomorphic. Thus $\Coker(A)\cong \Coker(\oplus A_i)$.
  Now invoke \S\ref{Sec.Background.Notations.Conventions}.v.
\epr

\beR
\bee[\bf i.]
\item
Assume $A\in \Mat{R}$ is homogeneous, i.e.,  all its entries are homogeneous, and of the same degree. Then we get a stronger statement:
 if $A\otimes R_\cm\sim \oplus \tA_i$ (for some $\{\tA_i\}$ over $R_\cm$) then
 $A\sim \oplus A_i$  by  $GL(m,R_0)\times GL(n,R_0)$.
\item
The locality of $R_0$ is important and the lemma does not hold if we localize just at the ideal $R_{>0}$. For example, suppose $R_0$ is a domain
 and there exists $A\in \Mats{R_0}$, with $det(A)\neq0$,  that is indecomposable.
 Then $(R_0)_{R_{>0}}$ is a field (in $R_{R_{>0}}$) and $A\otimes  R_{R_{>0}}\sim \one\otimes  R_{R_{>0}}$.

As another example, suppose $R_0$ is a   domain, and   $0\neq a_0\in R_0$, $0\neq a_1\in R_{>0}$ satisfy:
 $a_0\cdot a_1=0$. The $R_{>0}$-localization sends $a_1$  to $0\in R_{R_{>0}}$. Then, for any matrix $A\in \Mat{R}$,  the $R_{>0}$-localization of $a_1\cdot A$ is (trivially) diagonalizable.
\eee\eeR

\subsection{Decomposability over $\Spec(R)$ vs that over $\Proj(R)$ for graded rings}\label{Sec.Decompos.over.R.vs.over.Proj(R)}
 In this subsection  $R=\oplus_{d\in \N} R_d$, with $R_0$ a field. Assume the ideal $R_{>0}$ is generated (not necessarily finitely) by $R_1$.
 Take the  projective  scheme $\Proj(R)$, with the structure sheaf $\cO_{\Proj(R)}$, see pg. 76-77, 116-123 of \cite{Hartshorne}, or \S27.8-27.12 of \cite{Stacks}.

Recall some detail. The (not necessarily closed) points of $Proj(R)$ correspond to non-maximal homogeneous primes, $\cp\ssetneq R_{>0}$.
 The stalks of the structure sheaf at such points can be presented explicitly, $\cO_{(Proj(R),\cp)}=R_{(\cp)}\sset R_\cp$, where $R_{(\cp)}$ consists of
 the (equivalence classes of) fractions $\frac{a}{b}.$ Here $a\in R$, $b\in R\smin \cp$ are homogeneous and satisfy $deg(a)=deg(b)$. The basic
   affine charts are $\cU_g:=Proj(R)\smin V(g)$,
   for non-nilpotent elements $g\in R_1$. Their structure rings are $\cO_{Proj(R)}(\cU_g)=R_{(g)}\sset R[\frac{1}{g}]$,
    where $R_{(g)}:=\{\sum\frac{a_d}{g^d}|\ deg(a_d)=deg(g^d)\}$.

The correspondence $mod_{gr}R\rightsquigarrow Coh(Proj(R))$ associates to the graded module $M:=\Coker(A)$
  the coherent sheaf $\tM:=\widetilde{\Coker(A)}\in Coh(Proj(R))$.
 For a  matrix  $A=\sum A_d\in \Mat{\oplus R_d},$ graded in the sense of \S\ref{Sec.Background.Notations.Conventions}.v., and an   affine chart $\cU_g\sset Proj(R)$ we get the presentation
\beq\label{Eq.Presentation.Sheaf.Coker(A)}
\cO_{Proj(R)}(\cU_g)^n\stackrel{A_{(g)}}{\to}\cO_{Proj(R)}(\cU_g)^m\to \tM(\cU_g)\to0,\quad \text{where}\quad A_{(g)}=\sum \frac{A_d}{g^d}.
\eeq
Localizing this sequence at a (closed) point $\cp\in \cU_g$ we get the presentation of the stalk of the sheaf:
\beq\label{Eq.Sheaf.corresp.to.module.M}
\cO_{(Proj(R),\cp)} ^n\stackrel{A^{(\cp)}}{\to}\cO_{(Proj(R),\cp)} ^m\to \tM_\cp\to0.
\eeq

Here $A^{(\cp)}$ was obtained through the choice of the embedding $\cp\in \cU_g\sset Proj(R)$, but the morphism $A^{(\cp)}$ depends on $A,\cp$
 only.

 As $A$ is graded, all its determinantal ideals are homogeneous. Therefore $I_\bullet(A_{(g)})=I_\bullet(A)_{(g)}.$
 If $I_m(A)\!=\!J_1\cdot J_2$ then $I_m(A_{(g)})\!=\!J_{1,(g)} \cdot  J_{2,(g)}\!\sset\! R_{(g)}$ and
  $I_m(A^{(\cp)})=J^{(\cp)}_1\cdot J^{(\cp)}_2\sset \cO_{(Proj(R),\cp)}=R_{(\cp)}.$

\bthe\label{Thm.Decompos.Proj(R).vs.Local}
  Suppose  $A\in \Mat{R}$ is graded (\S\ref{Sec.Background.Notations.Conventions}.v)
  and   $I_m(A)=J_1\cdot J_2=J_1\cap J_2$   for some   homogeneous ideals.

\bee[\bf 1.]
\item If the conditions \eqref{Eq.Assumption.on.A} hold for $A\!\in\!\! \Mat{R}$ then they hold for all the stalks $A^{(\cp)}\!\!\in\!\! \Mat{\cO_{(Proj(R),\cp)}}.$
\item
  Suppose $H^0_{J_1+J_2}(M)=0$.
Assume there exists $x\in R_{>0}$ that is $M$-regular and  the base change $R\to \quots{R}{(x)}$ satisfies:
$I_m(\bA)=\bJ_1\cdot \bJ_2=\bJ_1\cap \bJ_2$ and $H^0_{\bJ_1+\bJ_2}(\bbM)=0$.
%\bei
%\item the natural embeddings $\quots{R}{(x)}\cdot H^0_{J_1+J_2}(M)\hookrightarrow H^0_{\bJ_1+\bJ_2}(\bbM)$ and
% $\quots{R}{(x)}\cdot H^0_{J_j}(\quots{R}{J_i}\otimes M)\stackrel{\de_2}{\hookrightarrow} H^0_{\bJ_j}(\quots{R}{J_i+(x)}\otimes M)$
%  of proposition \ref{Thm.Obstruction.Functoriality} are isomorphisms.
%\eei
\noindent Then  $A$ is $(J_1,J_2)$-decomposable  \iff for  each (closed) point $\cp\in   \P V(J_1)\cap \P V(J_2)\sset \Proj(R)$   the stalk
  $A^{(\cp)}$ is $ (J^{(\cp)}_1,J^{(\cp)}_2)$-decomposable.
\eee
\ethe
  Thus the graded decomposability problem in dimension $dim(R)$ is reduced to (many) local decomposability problems in dimension $dim(R)-1$.

 One can restate the conclusion for modules: ``A graded $R$-module $M$ is decomposable \iff all the stalks of the corresponding sheaf $\tM$ are decomposable".
\bpr
\bee[\bf 1.]
\item  The Fitting ideal sheaves of $\tM$ are the sheaves associated to the Fitting ideals of $M$.
 Indeed by \eqref{Eq.Presentation.Sheaf.Coker(A)} we identify the ideals on the basic affine opens:
\beq
Fitt_\bullet(\tM(\cU_g))=Fitt_\bullet(Coker(A_{(g)}))=Fitt_\bullet(M)_{(g)}=\widetilde{Fitt_\bullet(M)}(\cU_{g})\sset \cO_{Proj(R)}(\cU_g)=R_{(g)}.
\eeq
This identification of ideals is compatible with all the restrictions onto open subsets of $\cU_{g}$.
 Hence the identification of the sheaves of ideals,
 $Fitt_\bullet(\tM)= \widetilde{Fitt_\bullet(M)}\sset \cO_{Proj(R)}$.
 In particular we get: $Fitt_0(\tM)=\tJ_1\cdot \tJ_2$.

Let $f_i\in J_i$ be regular homogeneous elements of \eqref{Eq.Assumption.on.A}. We prove: their images in the stalks $ J^{(\cp)}_i\sset R_{(\cp)}$ are regular.
 It is enough to verify the regularity on the basic affine charts, as the localization $R_{(g)}\to R_{(\cp)}$ is exact.
 Take a chart, $\cp\in \cU_g\sset Proj(R)$, for some $g\in R_1$.
   Suppose $\frac{f_i}{g^{d_i}}\cdot \sum^N\frac{c_l}{g^l}=0\in R_{(g)}=\cO_{Proj(R)}(\cU_g)$. Then
    $ f_i \cdot \sum^N_l c_l g^{N-l}=0\in R[\frac{1}{g}]$. But $f_i\in R$ is regular in $R[\frac{1}{g}]$, therefore $\sum^N_l c_l g^{N-l}=0\in R[\frac{1}{g}]$.
     And thus $\sum^N\frac{c_l}{g^l}=0\in R_{(g)}$.

Finally we claim: $\tJ_1\cap \tJ_2=\tJ_1\cdot \tJ_2\sset \cO_{Proj(R)}$, i.e.,  these  sheaves of ideals are co-regular.
Namely, for each stalk: $J^{(\cp)}_1\cap J^{(\cp)}_2=J^{(\cp)}_1\cdot J^{(\cp)}_2\sset R_{(\cp)}$.
 By lemma \ref{Thm.coprime.ideals} it is enough to verify (locally):
 \beq
 Tor^{R_{(\cp)}}_1(\quots{R_{(\cp)}}{J^{(\cp)}_1},\quots{R_{(\cp)}}{J^{(\cp)}_2})=0 \quad
 \text{  for   all the points } \quad \cp\in \P V(J_1)\cap \P V(J_2).
  \eeq
    As the localization is exact, it is enough to verify this vanishing on the basic opens, i.e.
 $Tor^{R_{(g)}}_1(\quots{R_{(g)}}{J_{1,(g)}},\quots{R_{(g)}}{J_{2,(g)}})=0.$
   But this is the degree=0 part of  $Tor^{R[\frac{1}{g}]}_1(\quots{R[\frac{1}{g}]}{J_1[\frac{1}{g}]},\quots{R[\frac{1}{g}]}{J_2[\frac{1}{g}]})$.
 And the morphism $R\to R[\frac{1}{g}]$ is flat, therefore:
\beq
Tor^{R[\frac{1}{g}]}_1(\quots{R[\frac{1}{g}]}{J_1[\frac{1}{g}]},\quots{R[\frac{1}{g}]}{J_2[\frac{1}{g}]})=
 R[\frac{1}{g}]\otimes_R Tor^R_1(\quots{R}{J_1},\quots{R}{J_2})=0.
\eeq
Thus conditions \eqref{Eq.Assumption.on.A} hold for all the stalks $A^{(\cp)}.$
\item
 The part $\Rrightarrow$ is trivial. We prove the part $\Lleftarrow$.
\\
  By proposition \ref{Thm.Obstruction.to.Decomposability} it is enough to verify: $Q\!=\!0.$
 First we prove:  this obstruction module is supported on $V(R_{>0})\!\sset\! \Spec(R)$ only, i.e. just at one point.
  Then we prove: $depth(Q)\!\!>\!\!0.$ Together this implies  $Q\!=\!0.$

\bee[\hspace{-0.7cm}\bf Step 1.\!\!]
\item
 For each (closed) point $\cp\in Proj(R)$ take the stalk $M_{\cp}$ of \eqref{Eq.Sheaf.corresp.to.module.M}.
  By Step 1: $Fitt_0(\tM_\cp)=\tJ^{(p)}_1\cap \tJ^{(p)}_2=\tJ^{(p)}_1\cdot \tJ^{(p)}_2$.

In our case the exact sequence \eqref{Eq.Exact.Sequence.For.Decomposability} is  graded.
The functor $mod_{gr}R\to Coh(Proj(R))$ is exact, see e.g. \S27.8-27.12 of \cite{Stacks}. Therefore
 we get the sheaf analog of \eqref{Eq.Exact.Sequence.For.Decomposability}
  on $Proj(R)$. Localize it at $\cp$ to get the exact sequence of stalks:
\beq
0\to \tM_\cp \to     (\tM_1)_\cp  \oplus  (\tM_2)_\cp \to \tQ_\cp \to 0.
\eeq
Here $ (\tM_i)_\cp=(\widetilde{\quots{M'_i}{H^0_{J_j}(M'_i)}})_{(\cp)}=\quots{\tM'_{i,\cp}}{H^0_{\tJ^{(\cp)}_j}(\tM'_{i,\cp})}=(\tM_{\cp})_i$.
 Therefore, as $\tM_{\cp}$ is decomposable, proposition \ref{Thm.Obstruction.Functoriality} gives: $\tQ_{\cp}=0$.

 Thus the sheaf $\tQ\in Coh(Proj(R))$ vanishes. Therefore $Supp(Q)\sseteq V(R_{>0})\sset \Spec(R),$ i.e. the module $Q$ is supported at one point only.  In particular, either $depth(Q)=0$ or $Q=0$.

\item Take $x\in R$ as in the assumptions. As in the proof of proposition \ref{Thm.Obstruction.Functoriality}, we can either
 restrict equation \eqref{Eq.Exact.Sequence.For.Decomposability} to the hypersurface $V(x)\sset \Spec(R)$ (by applying $\quots{R}{(x)}\otimes$)
  or can write \eqref{Eq.Exact.Sequence.For.Decomposability} directly for $\bbM$.
 Then the diagram \eqref{Eq.Obstruction.Base.Change.Diagram} becomes:

\beq\label{Eq.diagram.on.Proj-R}
\bM
\oplusl_i Tor_1(\quots{R}{(x)},M_i)\!\to\! Tor_1(\quots{R}{(x)},Q)\to\!\!&\!\!\quots{R}{(x)}\otimes M\!\!&\!\!\stackrel{\pi}{\to}\!\!&\!\!
 \oplusl_i (\quots{R}{(x)}\otimes M_i)\!\!&\!\!\to\!\!&\!\!  \quots{R}{(x)}\otimes Q\!\!&\!\!\to\!\!&\!\!0
\\
&\rotatebox[]{-90}{$\isom{}$}&&\de_0\rotatebox[]{-90}{$\twoheadrightarrow$}&&\phi\rotatebox[]{-90}{$\twoheadrightarrow$}
\\
\hfill 0\!\!\to\!\!& \!\!\bbM\!\!&\!\!\stackrel{\bar\pi}{\to}\!\!&\!\!\oplus \bbM_i\!\!&\!\!\to\!\! &\!\!\bQ\!\!&\!\!\to\!\!&\!\!0
\eM
\eeq
(Note that $H^0_{J_1+J_2}(M)=0$ and  $H^0_{\bJ_1+\bJ_2}(\bbM)=0$.) \quad
 Then $\pi$ is injective, as $\bar\pi$ is injective.

We claim: $x$ is $M_i$-regular, i.e.,  $Tor_1(\quots{R}{(x)},M_i)=0$. Indeed, start from the presentation \eqref{Eq.Lifted.M.killed.torsion}.
 Suppose $x\cdot [\xi ]=0\in M_i$ for some $\xi\in M$. Then $ x\cdot \xi\in H^0_{J_j}(\quots{M}{J_i M})$. Therefore $J^d_j\cdot x\cdot \xi=0\in \quots{M}{J_i M}$ for some $d\gg1$.
  And then $J^{d+1}_j\cdot x\cdot \xi=0\in M$. But $x$ is $M$-regular, thus $J^{d+1}_j\cdot \xi=0\in M$. Hence $\xi\in H^0_{J_j}(M'_i)$,
   i.e.,  $[\xi ]=0\in M_i$.

Altogether we have: $Tor_1(\quots{R}{(x)},M_i)\!=\!0$ and  $\pi$ is injective. The diagram \eqref{Eq.diagram.on.Proj-R} gives: $Tor_1(\quots{R}{(x)}\!,\!Q)\!=\!0$.
 Thus if $Q\!\neq\!0$ then $x$ is $Q$-regular,
 in contradiction with $depth( Q)\!=\!0.$
\eee
\mbox{\!Altogether $Q\!=\!0$. \!By \!proposition \ref{Thm.Obstruction.Functoriality}:  $\Coker(A)\!=\!\oplus\! M_i$.
 \!Thus (\S\ref{Sec.Background.Notations.Conventions}.v) $A$ is \!decomposable \!over $R.$
\epr}\eee

 \beR\label{Rem.3.13}
 The assumption $\bJ_1\cdot\bJ_2=\bJ_1\cap \bJ_2$ implies $\empty\neq \P V(J_1)\cap \P V(J_2)\sset \Proj(R).$
   Indeed, if $\P V(J_1)\cap \P V(J_2)=\empty,$ then $\sqrt{J_1+J_2}=R_{>0}.$ But then the intersection $V(\bJ_1)\cap V(\bJ_2)$ cannot be proper.

 The proposition does not hold when $\empty\!=\! \P V(J_1)\!\cap\! \P V(J_2)\!\sset\! \Proj(R)$. For example, let $R\!=\!\k[x,y],$
 for a field
$\k\!=\!\bar\k.$
 Take $A\!\in\! \operatorname{Mat}_{2\times 2}(R_d)$, whose elements are generic homogeneous polynomials of degree $d$.
 Then $\det(A)$ splits into co-prime linear factors.
 The stalks $\{A^{(\cp)}\}$ are trivially decomposable.
 But for $d\!\ge\!3$ the matrix is indecomposable. (Note that $I_1(A)$ is minimally generated by 4 elements.)
Compare this to the   decomposition of modules with non-connected support, \S\ref{Sec.Background.Decomposability.for.non.Connected.Support}.
  \eeR

The technical assumptions of part 2 of theorem \ref{Thm.Decompos.Proj(R).vs.Local} become simple in the case of square matrices.
\bcor\label{Thm.Decompos.Proj.vs.Local.square.case}
Let $det(A)\!=\!f_1\!\cdot\! f_2$, regular and coprime. Suppose $depth(\quots{R}{(f_1,f_2)})\!>\!0.$
  Then $M$ is $f_1,\!f_2$-decomposable \iff $\tM_\cp$ is $f_1^{(\cp)}\!, f_2^{(\cp)}$-decomposable
 at each (closed) point  $\cp\!\in\! \P V(f_1)\!\cap\! \P V(f_2)\!\sset\! Proj(R).$
\ecor
\bpr By lemma \ref{Thm.coprime.ideals} we have $H^0_{(f_1,f_2)}(M)=0.$
Take $x\in R$ that is $\quots{R}{(f_1,f_2)}$-regular. We can assume that $x$ is homogeneous. Then $f_1,f_2,x$ is a regular sequence. As $R$ is graded,
 and $f_1,f_2,x$ are homogeneous, the sequence $x,f_1,f_2$ is regular as well. Therefore $\bar{f}_1,\bar{f}_2\in \quots{R}{(x)}$ are coprime.
 Therefore  and   $H^0_{(\bar{f}_1,\bar{f}_2)}(\bbM)=0$, by lemma \ref{Thm.coprime.ideals}.

Moreover, $x$ is $M$-regular. Indeed, if $x\cdot \xi\in Im(A),$ then $x\cdot Adj(A)\xi\in (f_1f_2)\cdot R^n.$ Thus  $  Adj(A)\xi\in (f_1f_2)\cdot R^n$, hence $\xi\in Im(A)$.

 Now apply theorem \ref{Thm.Decompos.Proj(R).vs.Local}.
\epr

\beR
Recall Max Noether's fundamental theorem. Given homogeneous polynomials $f,g,h\in \k[x_0,x_1,x_2]$, $\k=\bar\k$, with $f,g$ coprime.
 Then $h\in (f,g)\sset \k[x_0,x_1,x_2]$ \iff this
 holds for the localizations,  $h_\cp\in (f_\cp,g_\cp),$ at all the points of the intersection $V(f)\cap V(g)\sset \P^2$.
   Theorem \ref{Thm.Decompos.Proj(R).vs.Local} is the natural analogue of this theorem for the decomposability question of matrices over graded rings.
\eeR

\section{Decomposability criteria for square matrices}\label{Sec.Decomposability}
 \bthe\label{Thm.Decomposition.Square.Matrices}
  Let $R$ be a commutative, unital ring and let $A\in \Mats{R}$.   Suppose   $\det(A)=f_1 \cdots f_r\in R$, where
    $\{f_i\}$ are  regular and pairwise co-prime, see \S\ref{Sec.Background.Coprime.elements.Coregular.ideals}.

\bee[1.]
\item     $A$ is stably-$(f_1, \dots,f_r)$-decomposable \iff $I_{n-1}(A)\sseteq \sum^r_{j=1} (\prod_{i\neq j} f_i)$.
\item Assume   $R$ is a local ring, or $R$ is graded and $A$ is graded. Then
 $A$ is $(f_1, \dots,f_r)$-decomposable \iff $I_{n-1}(A)\sseteq \sum^r_{j=1} (\prod_{i\neq j} f_i)$.
\eee
 \ethe
For another presentation of the condition $I_{n-1}(A)\sseteq \sum^r_{j=1} (\prod_{i\neq j} f_i)$  see lemma \ref{Thm.Coprime.Elements}.

\

This theorem is proved in \S\ref{Sec.Proof.Decomp.Linear.Algebra}, using auxiliary results of \S\ref{Sec.Alomost.Projectors}.
 Then we give examples, \S\ref{Sec.Decomposition.Examples}, and discuss the (necessity of the) assumptions, \S\ref{Sec.Decomposition.Remarks}.
 Then come the first applications to the diagonal reduction of linear determinantal representations/tuples of matrices.

\subsection{Block-diagonalization of (``almost"-)projectors  over a ring}\label{Sec.Alomost.Projectors}
In this subsection we impose the condition on $R$:
\beq\label{Eq.Fin.Gen.proj.are.free}
\text{any finitely-generated projective module over $R$ is free.}
\eeq
It holds, e.g., for local rings and for $S[x]$,
 with $S$ local Noetherian, $dim(S)\le2$, see \cite[\S A.3.2]{Eisenbud-book}.
\bel\label{Thm.Lemma.on.Projectors}
Suppose $R$ satisfies \eqref{Eq.Fin.Gen.proj.are.free}, while the operators $P_1,\dots,P_r \in \Mats{R} $ satisfy: $\sum P_i=\one_{n\times n}$ and $P_i P_j=\zero_{n\times n}$
 for $i<j$.
 Then $\{P_i\}_i$ are simultaneously diagonalizable. Namely, there exists $U\in GL(n, R)$ satisfying:
 $\sum x_i U P_i U^{-1}=\oplus x_i U P_i U^{-1}=\oplus x_i \one_{n_i}$. (Here $\{x_i\}$ are indeterminates.)
\eel
\bpr
First we establish the case $r=2$.
\bee[\bf Step 1.]
\item We verify the standard projector properties: $P^2_i=P_i$ and $P_iP_j=\zero$ for $i\neq j$. Indeed:
\beq
P^2_1=P_1(P_1+P_2)=P_1\cdot\one=P_1,\quad\quad\quad P^2_2=(P_1+P_2)P_2=\one\cdot P_2=P_2
\eeq
From here one gets: $P_1=(P_1+P_2)P_1=P^2_1+P_2P_1=P_1+P_2P_1$, hence $P_2P_1=\zero$.
\item
Consider $P_1,P_2$ as endomorphisms of the free module $F= R^n$. Define $F_i=P_i(F)$.
 We claim: $ F_1\oplus F_2=F$. Indeed, by the definition: $P_i(F_j)=P_i(P_j F)=0$ for  $i\neq j$.
 Thus: $P_i(F_1\cap F_2)=0$ for any $i$. But then: $\one (F_1\cap F_2)=(P_1+P_2)(F_1\cap F_2)=0$.

 Besides: $F_1\oplus F_2= P_1(F)\oplus P_2(F)=\one(F)=F$.  Therefore $F=F_1 \oplus F_2.$
  Hence $F_1,F_2$ are projective submodules of $F$.
  But then, by the initial assumption on $ R$,  they are free.

\medskip

Finally, take some bases $ v^{(j)}_\bullet$ of $F_j$, for $j=1,2$.
 The change  from the standard basis of $F$ to  the basis  $ v^{(1)}_\bullet, v^{(2)}_\bullet$  gives the needed
   diagonalizing transformation $P_i\to UP_iU^{-1}$.
\eee

For the case $r>2$ we apply the ($r=2$)-argument to the operators $P_1,\sum_{i=2}^r P_i$. Then we restrict to the submodule $\sum_{i=2}^r P_i(F)$ and iterate.
\epr
\beR
 The assumption \eqref{Eq.Fin.Gen.proj.are.free} is necessary. Suppose  there exists a (f.g.)
  projective but non-free module $F_1\in\mod R$. Complement it to a free module,
 $F_1\oplus F_2=R^n$. Define the projection homomorphisms $P_i:R^n\to F_i\subset R^n$ by
 $R^n\ni s=s_1+s_2\to s_i\in F_i$. (Here the decomposition $s=s_1+s_2$ is unique.)
Then $P_1P_2=\zero$ and $P_1+P_2=\one$. But $\{P_i\}$ cannot be brought to the prescribed form. (This would imply the freeness of $F_i$.)

As an explicit example take $R=\quots{\k[x]}{(x^2-1)}$, $char(\k)\neq2$, and $P_\pm:=(\frac{1\mp x}{2})\oplus(\frac{1\pm x}{2})\in \operatorname{Mat}_{2\times 2}(R)$.
\eeR

Lemma \ref{Thm.Lemma.on.Projectors} dealt with projectors. Now we treat the ``almost-projectors".
  Let $\ca\ssetneq R$ be an  ideal in a local ring, and suppose $(R,\ca)$ is a Henselian pair, see \cite[\S15.11]{Stacks}.
\bel\label{Thm.Lemma.on.almost.projectors}
 Suppose $P_1,\dots,P_r\!\in\! \Mats{R} $ satisfy:
 $\sum P_i\!=\!\one$ and $P_iP_j\!\in\! \Mats{\ca}$ for $i\!<\!j$.
 Then $\{P_i\}_i$ are simultaneously block-diagonalizable. Namely  there exists $U\in GL(n, R)$ satisfying
  \beq\label{Eq.Almost.Projectors}
  \sum x_i U P_iU^{-1}=\oplus_i \big(x_i \one_{n_i}+\sum_j x_j A_{ij}\big).
  \eeq
Here $\{x_i\}$ are indeterminates,   $A_{ij}\in \operatorname{Mat}_{n_i\times n_i}(\ca)$ and $\sum_j  A_{ij}=\zero$ for each  $i$.
\eel
%  \square \!\!\! {\substack{\! \ca\\}} \oplus \square \oplus \square \!\!\! \substack{\! \ca\\}$.
% \[U P_1U^{-1}=\bpm \square&\zero\\\zero&\square_\ca\epm\quad and\quad U P_2U^{-1}=\bpm \square_\ca&\zero\\\zero&\square\epm,\]
% Here all the entries in the blocks $\square \!\!\! \substack{\! \ca\\}$ belong to $\ca$ and the blocks are of the corresponding sizes.
\bpr
 First we establish the case $r=2$.
Take the quotient $R\stackrel{\phi}{\to} \quots{R}{\ca}$. We get the projectors: $\phi(P_1)+\phi(P_2)=\one$, $\phi(P_1)\cdot \phi(P_2)=\zero$.
 Thus we can assume
$\phi(P_1)=\bbm \one&\zero\\\zero&\zero\ebm$ and  $\phi(P_2)=\bbm \zero&\zero\\\zero&\one\ebm$,
by lemma \ref{Thm.Lemma.on.Projectors}.
 (Note that the ring $\quots{R}{\ca}$ is local, thus the condition  \eqref{Eq.Fin.Gen.proj.are.free} holds.)

Now we eliminate the off-diagonal blocks of $P_i$.
 For $P_1$ it is enough to resolve  the condition
\beq\label{Eq.Almost.Projectors.eq.to.resolve}
(\one+\tU)P_1(\one+\tU)^{-1}=\bbm *&\zero\\\zero&*\ebm\quad {\rm for} \quad \tU\in \Mats{\ca}.
\eeq

 Present $P_1$ and $\tU$ in the block-form, then condition \eqref{Eq.Almost.Projectors.eq.to.resolve} will follow from
\beq
\bbm \one&\tU_{12}\\\tU_{21}&\one\ebm\bbm \one+P_{11}&P_{12}\\P_{21}&P_{22}\ebm=\bbm D_1&\zero\\\zero&D_2\ebm\bbm \one&\tU_{12}\\\tU_{21}&\one\ebm.
\eeq
Here  all the entries of $\{P_{ij}\}$, $\tU_{12}$, $\tU_{21}$ belong to $\ca.$
 (And $\tU_{12},\tU_{21},D_1,D_2$ are unknowns.)
 Equality of the diagonal blocks gives $D_1=\one+P_{11}+\tU_{12}P_{21}$ and $D_2=\tU_{21}P_{12}+P_{22}$. Substitute these $D_1,D_2$ into the off-diagonal equations, then we have to resolve:
 \beq\label{Eq.equation.to.resolve}
 P_{12}+\tU_{12}P_{22}=(\one+P_{11}+\tU_{12}P_{21})\tU_{12},\quad\quad
   \tU_{21}(\one+P_{11})+P_{21}=(\tU_{21}P_{12}+P_{22})\tU_{21}.
 \eeq
This is a system of polynomial equations in $\tU_{12},\tU_{21}.$
 It is  solvable (by the implicit function theorem), as $(R,\ca)$ is a henselian pair.  Indeed, e.g., for the map
 $F(X)= P_{12}+XP_{22}-(\one+P_{11}+XP_{21})X$ the derivative $\frac{\di F}{\di X}|_{X=\zero}$
 is invertible, being of the form $\one+(\ca)$.

Thus $P_1$ is brought by conjugation to the form $D_1\oplus D_2$, where $D_1\in \one\oplus \operatorname{Mat}_{n_1\times n_1}(\ca)$,
  $D_2\in  \operatorname{Mat}_{n_2\times n_2}(\ca)$.
 The prescribed block-diagonal form of $P_2$   follows, as the conjugation preserves the conditions  $P_1+P_2=\one$, $P_1P_2\in \Mats{\ca}$.
  Thus we have reached the form  of \eqref{Eq.Almost.Projectors}.

\

For the case $r>2$  apply the ($r=2$) argument to $P_1,\sum^r_{i=2}P_i$, and iterate.
\epr
\beR
The assumption ``the pair $(R,\ca)$ is local Henselian" cannot be weakened to ``$ R $ is a local ring".
 For example, consider the $2\times 2$ case:
 \beq
 P_1=\bbm 1-f&g\\g&f\ebm,
  \quad\quad P_2=\bbm f&-g\\-g&1-f\ebm,
   \quad \quad  f,g\in(\ux)\sset \k[\ux]_{(\ux)},\quad char(\k)\neq2.
\eeq
 Then the block-diagonalization means:
\beq
P_1\rightsquigarrow \bbm \frac{1+\sqrt{1+4(g^2+f^2-f)}}{2}&0\\0&\frac{1-\sqrt{1+4(g^2+f^2-f)}}{2}\ebm.
\eeq
The resulting matrix is over   $\k\bl\ux\br$ but not over $\k[\ux]_{(\ux)}$, regardless of how large are the orders of $f,g$.
\eeR

\subsection{The proof of theorem \ref{Thm.Decomposition.Square.Matrices}}\label{Sec.Proof.Decomp.Linear.Algebra}

The directions $\Rrightarrow$ are trivial. We prove the directions $\Lleftarrow$.

It is enough to prove part 1. Then part 2  follows  by \S\ref{Sec.Background.Notations.Conventions}.v.

By corollary \ref{Thm.reduction.to.complete.local.ring.1} and part 2.iii. of lemma \ref{Thm.coprime.ideals} we can assume: $(R,\cm)$
 is local and Henselian. Then (by \S\ref{Sec.Background.Notations.Conventions}.v) we can pass to the minimal presentation, assuming $A\in \Mats{\cm}$.

\bee[\!\!\!\bf Step 1.\!]
\item

By the assumption, $ Adj(A)=\sum_j (\prod_{i\neq j}f_i)B_i$, for some matrices  $B_i\in\Mats{R}$.
This decomposition is not unique  due to the freedom
\beq\label{Eq.Freedom.Shuffle}
B_i\to B_i+f_i Z_i\quad \text{ for }\quad \sum Z_i=\zero.
\eeq
 We will
use this freedom later.

From this presentation of $Adj(A)$ we get:
  $\det(A)\cdot \one=A\cdot Adj(A)=\sum_j (\prod_{i\neq j}f_i)AB_i.$
 As $\{f_i\}$ are regular and co-prime, we have  $(\prod_{i\neq j}f_i)\cap (f_j)=(\prod^r_{i=1}f_i)$, lemma \ref{Thm.Coprime.Elements}.
 Therefore we get $(\prod_{i\neq j}f_i)AB_i\in (\det(A))\cdot \Mats{R}$, and thus $AB_i\in (f_i)\cdot \Mats{R}$.
Therefore we define the matrices $\{P_i\}$, $\{Q_i\}$ by $f_iP_i:=AB_i$ and   $f_iQ_i:= B_iA$.
By their definition:
$\sum P_i=\one$ and $\sum Q_i=\one$.

 We want to conclude:  $\oplus P_i=\one$ and $ \oplus Q_i=\one$.
The key ingredient is the identity:
\beq
f_i B_j P_i= B_jA B_i=f_jQ_j B_i,\quad\quad\text{for} \quad \quad  i\neq j.
\eeq
 As $f_i,f_j$ are regular and co-prime,    $B_jP_i$ must be divisible by $f_j$, i.e.,  $\frac{B_j P_i}{f_j}\!=:\!\tZ_{ji}\!\in\! \Mats{R}.$
   Then, for any $i\!\neq\! j$ we get:
$P_jP_i\!=\!\frac{AB_j P_i}{f_j}\!=\!A\tZ_{ji}$.
Similarly, $Q_jQ_i\!=\!\tZ_{ji}A$ for
 any $i\!\neq\! j.$
Thus $\{P_i\}$ and $\{Q_i\}$ are  almost projectors in the sense of lemma \ref{Thm.Lemma.on.almost.projectors}, for $\ca\!=\!\cm.$

\item
  \mbox{Assume $r\!=\!2$.
The transformation $A\!\to\! UAV$ results in $B_i\!\to\! Adj(V)\!\cdot\! B_i\!\cdot \!Adj(U).$ This implies:}
\beq
(\frac{AB_i}{f_i}=)\quad\quad P_i\to \det(UV)\cdot U P_i U^{-1} \quad\quad \text{ and } \quad\quad  (\frac{B_iA}{f_i}=)\quad\quad Q_j\to \det(UV)\cdot V^{-1} Q_j V.
\eeq
The conditions $P_1P_2=A\tZ$ and $Q_1Q_2=\tZ A$ are transformed into:
\beq
P_1P_2=\frac{1}{det(UV)^{2}}UA\tZ U^{-1},\quad\quad\quad
Q_1Q_2=\frac{1}{det(UV)^{2}}V\tZ AV^{-1}.
\eeq

Thus, using lemma \ref{Thm.Lemma.on.almost.projectors}, we can assume (using $U$):
\beq
P_1=\bbm \one-Z_1&\zero\\\zero&Z_2\ebm,\quad\quad\quad P_2=\bbm Z_1&\zero\\\zero&\one-Z_2\ebm.
\hspace{2cm}  \text{Here $Z_i\in \operatorname{Mat}_{n_i\times n_i}(\cm)$.}
\eeq
Similarly (by using $V$) we have the block-diagonalization of  $Q_1,Q_2$.

The condition $P_1P_2=A\tZ$ gives:
\beq
\bbm (\one-Z_1)Z_1&\zero\\\zero&Z_2(\one-Z_2)\ebm\!=\!A\tZ,  \text{ and thus }
\bbm Z_1&\zero\\\zero&-Z_2\ebm\!=\!A\tZ\!\bbm(\one-Z_1)^{-1}&\zero\\\zero&-(\one-Z_2)^{-1}\ebm.
\eeq
Now apply  the freedom of equation \eqref{Eq.Freedom.Shuffle},
$B_1\to B_1+f_1 Z$, $B_2\to B_2-f_2 Z.$ This amounts to: $P_1\to P_1+AZ$ and
$P_2\to P_2-AZ$. Thus  we choose
\beq
Z=\tZ\bbm(\one-Z_1)^{-1}&\zero\\\zero&-(\one-Z_2)^{-1}\ebm \quad \text{ to get:}
 \quad \quad   P_1\to\bbm \one&\zero\\\zero&\zero\ebm,\quad P_2\to\bbm \zero&\zero\\\zero&\one\ebm.
\eeq
 Take  $P_1,P_2$ in this form,  thus  $\zero=P_1P_2=\frac{AB_1AB_2}{f_1f_2}$. Therefore   $B_1AB_2=\zero$ and   $Q_1Q_2=\zero$.

By lemma \ref{Thm.Lemma.on.Projectors} we can assume $Q_1=\one\oplus \zero$ and $Q_2=\zero\oplus \one$.

Finally, use the original definition of $P_i$ and $Q_i$, to write: $\frac{f_i}{f} Adj(A)P_i\! =\!B_i\!=\!\frac{f_i}{f}Q_i Adj(A)$.
This gives $Adj(A)\!=\!f_2 B_1\!\oplus\! f_1 B_2$. Therefore $A\!=\!A_1\!\oplus\! A_2$, with $Adj(A_i)\!=\!B_i,$ and $det A_i\!=\!f_i.$

\

For $r>2$, apply the previous argument to $P_1,\sum_{i=2}^r P_i$. Decompose $A=A_1\oplus \tA_2$. Then restrict to the subspace corresponding to $\tA_2$
 and iterate.
\epr
\eee

\subsection{The first examples}\label{Sec.Decomposition.Examples}
\bex\label{Ex.Decomposability.Square.Matrices.Cases}
\bee[\bf i.]
\item  Suppose $R$ is a principal ideal domain, e.g., $R=\k[x],$ $\k[\![x]]\!$, for an arbitrary field $\k,$ or $\k\{x\}$ for $\k$ a normed field. Let $A\in \Mats{R}$,
  $\det(A)=f_1f_2$, with $f_1,f_2$ regular and co-prime.
 In this case $(f_1)+(f_2)=R$, as this ideal must be principal and cannot be proper (otherwise $ f_1,f_2$ are not co-prime). Thus
  the condition $I_{n-1}(A)\sseteq (f_1)+(f_2)$ is empty. Hence $A\stackrel{s}{\sim} A_1\oplus A_2$ with $\det(A_i)=f_i$. By iterating this procedure we get: if
   $\det(A)=\prod f^{p_i}_i$ is the decomposition into irreducible (non-invertible) pairwise co-prime elements then $A\stackrel{s}{\sim} \oplus A_i$, with $\det(A_i)=f^{p_i}_i$.

    This establishes   ``a half" of the classical Smith normal form for matrices over PID's.

\item  Let $A=\bbm y&x^k\\x^l&y\ebm$ with $R=S[[x,y]]$, $S$ being any (commutative, unital) ring.  Then $\det(A)=y^2-x^{k+l}$ is reducible \iff $k+l\in2\Z$.
 Suppose $k+l\in2\Z$, and $2\in R$ is invertible. Then $A$ is decomposable \iff $k=l$.  This goes by verifying the condition $I_1(A)\sseteq (y-x^{\frac{k+l}{2}})+(y+x^{\frac{k+l}{2}})$.

\item A bit more generally, let $R=S[[\ux]]$, with $\ux$ a multi-variable. Take a matrix of monomials,
  $A=\{a_{ij}\}=\{\ux^{\ud_{ij}}\}\in \operatorname{Mat}_{2\times 2}(R)$. Here $\{\ud_{ij}\} $ are  vectors of natural numbers.
    Assume all the coordinates of
   the vectors    $\ud_{11}+\ud_{22}$,    $\ud_{12}+\ud_{21}$ are even numbers, then
  \beq
  \det(A)=(\ux^\frac{\ud_{11}+\ud_{22}}{2}-\ux^\frac{\ud_{12}+\ud_{21}}{2})(\ux^\frac{\ud_{11}+\ud_{22}}{2}+\ux^\frac{\ud_{12}+\ud_{21}}{2})=:f_-\cdot f_+\in R.
\eeq
Assume $f_-,f_+$ are co-prime. We claim: $A$ is $(f_-,f_+)$-decomposable \iff \big($\ud_{11}=\ud_{22}$ and $\ud_{12}=\ud_{21}$\big).
 Indeed, $(f_-)+(f_+)=(\ux^\frac{\ud_{11}+\ud_{22}}{2})+(\ux^\frac{\ud_{12}+\ud_{21}}{2})$.
 Then the necessary condition for the decomposability, $I_1(A)\sseteq (f_-)+(f_+)$, reads:
\beq
\forall (i,j):\ \text{ either } \ud_{ij}\ge\frac{\ud_{11}+\ud_{22}}{2} \text{ or } \ud_{ij}\ge\frac{\ud_{12}+\ud_{21}}{2}.
\eeq
(Here $\ge$  means the  inequality in  each coordinate of these vectors.) Observe:
\bei
\item
 if  $\ud_{11},\ud_{22}\ge \frac{\ud_{12}+\ud_{21}}{2}$ then $f_-,f_+$ are not coprime;
 \item
  if  $\ud_{11}\ge \frac{\ud_{11}+\ud_{22}}{2}$ and $\ud_{22}\ge \frac{\ud_{12}+\ud_{21}}{2}$ then $\frac{\ud_{11}+\ud_{22}}{2}\ge \frac{\ud_{12}+\ud_{21}}{2}$, thus
  $f_-,f_+$ are not coprime;
\item the only remaining case is: $\ud_ {11},\ud_{22}\ge \frac{\ud_{11}+\ud_{22}}{2}$, implying $\ud_{11}=\ud_{22}$;
\item and similarly for $\ud_{12},\ud_{21}$.
\eei
Altogether we get:   $\ud_{11}=\ud_{22}$ and $\ud_{12}=\ud_{21}$.

Vice versa, if  $\ud_{11}=\ud_{22}$ and $\ud_{12}=\ud_{21}$ then $I_1(A)= (f_-)+(f_+)$.
 And thus $A$ is $(f_-,f_+)$-decomposable, by theorem \ref{Thm.Decomposition.Square.Matrices}.

\item Examples ii., iii. were for the ring of power series, $S[\![\ux]\!],$ though $A$ was a matrix of monomials. Invoking Lemma \ref{Thm.Decomposability.Graded-to-Local.reduction} we get in this case: $A$ is $(f_+,f_-)$-decomposable also polynomially, over $S[\ux].$

  \item Yet more generally, suppose   $A\in \operatorname{Mat}_{2\times 2}(R)$ satisfies:  $\det(A)=f_1f_2$, with $f_1,f_2$ regular and co-prime.
   Then $A$ is (stably-)$(f_1,f_2)$-decomposable
   \iff  the entries $a_{11},a_{12},a_{21},a_{22}$ all belong to the ideal $(f_1)+(f_2)\sset R$.
\eee
\eex

An immediate consequence of theorem \ref{Thm.Decomposition.Square.Matrices} is:
\bcor\label{Thm.Cor.4.7}
 Let $(R,\cm)$ be a local ring and $A\in \Mats{\cm}$. Suppose $\det(A)=f_1f_2$ with $f_1,f_2$   regular and co-prime, and $(f_1)+(f_2)\supseteq \cm^{n-1}$.
 Then $A$ is $(f_1,f_2)$-decomposable.
\ecor
We recall some cases when the condition $(f_1)+(f_2)\supseteq \cm^{n-1}$ holds.  As $A\in \Mats{\cm}$ the $\cm$-order of $\det(A)$ is at least $n$,
 assume $ord_\cm(\det(A))=n$.
  This corresponds to the {\em maximal} determinantal representations, or a representation of {\em maximal} corank,
  see \cite{Kerner-Vinnikov.Det.Reps.Global}.
  In this case the module $\Coker(A)$ is Ulrich-maximal, being minimally generated by $n$ elements, see  \cite{Ulrich84}.

 Assume  $ord_\cm f_i=n_i$, with   $n=n_1+n_2$.
 Take the leading terms  $l.t._{n_1}(f_1), l.t._{n_2}(f_2)\in gr_\cm(R),$ in the associated graded ring. By Nakayama it is enough to verify:
 \beq\label{Eq.to.verify}
 \quots{\cm^{n-1}}{\cm^n}\sseteq (l.t._{n_1}(f_1))+ (l.t._{n_2}(f_2))\sseteq gr_\cm(R).
 \eeq
Thus we study the $\quots{R}{\cm}$-vector subspace
\beq
 l.t._{n_1}(f_1)\cdot \quots{\cm^{n-1-n_1}}{\cm^{n-n_1}}+  l.t._{n_2}(f_2)\cdot\quots{\cm^{n-1-n_2}}{\cm^{n-n_2}}
\sseteq\quots{\cm^{n-1}}{\cm^n}.
\eeq
Assume $l.t._{n_1}(f_1), l.t._{n_2}(f_2)\in gr_\cm(R)$ is a regular sequence, then this vector subspace is the direct sum, of dimension
 $\dim \quots{\cm^{n-1-n_1}}{\cm^{n-n_1}}+\dim \quots{\cm^{n-1-n_2}}{\cm^{n-n_2}}$.

\bex\label{Ex.Decomposability.Square.Matrices.Maximal.Corank.origin}
\bee[\bf i.]
\item Let $(R,\cm)$ be a regular local Noetherian ring of Krull dimension 2. Then the Hilbert function $p_R(j):=dim\quots{\cm^{j-1}}{\cm^{j}}$
 satisfies   $p_R(n_1)+p_R(n_2)=p_R(n_1+n_2)$.
 For $\det(A)\!=\!f_1f_2$ we get two curve-germs, $V(f_1),V(f_2)\!\sset\! \Spec(R)$.
 The condition ``$jet_{n_1}(f_1), jet_{n_2}(f_2)\!\in\! gr_\cm(R)$  are co-prime" means: the  tangent cones of $V(f_1),V(f_2)$ have no common line.
 Suppose the $\cm$-order of $\det(A)$ equals $n=n_1+n_2$, then the condition \eqref{Eq.to.verify} holds.
  Then $A\sim A_1\oplus A_2$, with $\det(A_i)=f_i$.
\item More generally, suppose for $(R,\cm)$ the Hilbert function $p_R(j)$ satisfies   $p_R(n_1)+p_R(n_2)=p_R(n_1+n_2)$ for $n_1,n_2\gg1$.
 (A typical example is a two-dimensional local ring whose integral closure is regular.) Then the condition \eqref{Eq.to.verify}
   holds for $n_1,n_2\gg1$, when $jet_{n_1}(f_1), jet_{n_2}(f_2)$ are co-prime.
\eee
\eex
\beR
\bee[\bf i.]
\item
The assumption of maximal corank,  $ord_\cm(\det(A))=n$, is vital. For example, let $R=\k[[x]]$ and $A\in \operatorname{Mat}_{2\times2}(\cm^N)$ be a matrix
of homogeneous forms of degree $N$ in two variables and $\k=\bar\k$. Then $\det(A)$ necessarily splits. But if $N\ge3$ then usually
$I_1(A)$ cannot be generated by fewer than 4 elements. Hence $A$ is indecomposable, not even equivalent to an upper-block-triangular form.
\item
The regularity of $ R$ is important. If $dim_\k(\quots{\cm}{\cm^2})>2$ then usually $(f_1)+(f_2)\not\supseteq\cm^{n-1}$, even when the ideal $(f_1)+(f_2)$
 is $\cm$-primary.

\item The condition  $(f_1)+(f_2)\supseteq \cm^{n-1}$ does not hold when $dim(R)>2$, even for $n\gg1$.
\eee
\eeR

\subsection{Remarks on the assumptions of theorem \ref{Thm.Decomposition.Square.Matrices}}\label{Sec.Decomposition.Remarks}

\bee[\bf i.]

\item\label{Ex.Mutually.Prime.Essential}

 The condition ``$f_1,f_2$ are co-prime in $R$" is essential. It is the analogue of the condition on
distinct eigenvalues when diagonalizing a  matrix over a field.
 The theorem does not hold if $f_1,f_2$ are not co-prime.
As an example, take some matrix factorization, $A_1A_2=f\one$, for a non-zero divisor $f\in R$,
 such that $\det(A_i)=f^{p_i}$. Then $A:=\bbm A_1&B\\\zero&A_2\ebm$ satisfies: $\det(A)=f^{p_1+p_2}$ and
\beq
Adj(A)\!=\!\bbm \det(A_2)\!\cdot\! Adj(A_1)& -Adj(A_1)\!\cdot\! B\!\cdot \! Adj(A_2)\\\zero&\det(A_1)\!\cdot\! Adj(A_2)
\ebm\!=\!
\bbm f^{p_1+p_2-1}A_2&-f^{p_1+p_2-2}A_2BA_1\\\zero&f^{p_1+p_2-1}A_1\ebm\!.
\eeq
Thus, for $p_1,p_2\ge2$, we have the inclusion $I_{n-1}(A)\sset (f^{p_1+p_2-2})\sset(f^{p_1},f^{p_2})$.
But $A$ is not   equivalent to a block-diagonal matrix, as we made no assumptions on $B$.

 \item The condition $(f_1)\cap (f_2)=(f_1f_2)$ implies: the hypersurfaces $V(f_1),V(f_2)\sset \Spec(R)$ have no common component.
  Here $V(f_i)$ are taken as subschemes,  not just   the zero sets.
   And the condition is on  all the  components, including the embedded components.

   For example, for $R=\quots{\k[x,y,z]}{(z^2,z(x-y))}$ take $f_1=x$, $f_2=y$. The intersection of the {\em sets} $V(x)\cap V(y)\sset\k^2_{xy}$ is proper,
    but $(x)\cap (y)\neq (x\cdot y)$ inside $R$.

Another example (showing that all the closed points of a scheme should be checked)
  is    $R=\R[x,y]_{(x,y)}$ with $f_1=x(x^2+y^2)$, $f_2=y(x^2+y^2)$. Here $(f_1)\cap (f_2)\supsetneq(f_1f_2)$.

\item
 The condition  $I_{n-1}(A)\sseteq (f_1)+(f_2)$ is necessary for decomposability, by the direct check of $I_{n-1}(A_1\oplus A_2)$.
  The geometric meaning of this condition, when the subscheme $V(f_1,f_2)\sset \Spec(R)$ is reduced, is:
  $A$ is of $corank\ge 2$ on the locus $V(f_1)\cap V(f_2)\sset \Spec(R)$.

\item
We do not assume that $R$ is Noetherian. Thus theorem \ref{Thm.Decomposition.Square.Matrices} works e.g., for the rings of continuously-differentiable functions,
 $C^r(\cU)$, for $1\le r\le \infty$  and $\cU\sset \R^p$, or their germs along  closed subsets,  $C^r(\cU,Z)$, resp. their quotients by ideals,
   $\quots{C^r(\cU)}{I}$,  $\quots{C^r(\cU,Z)}{I}$.
   \bei
   \item Here the condition ``$f$ is regular" is easy to verify.
 Recall that   $f\in C^r(\cU)$ is not a zero-divisor \iff its zero locus, $V(f)\sset \cU$,  has empty interior.
 \item  But the co-primeness of $f_1,f_2$ is a restrictive condition. For example, let $R=C^\infty(\R^p,o)$ and $f_1,f_2\in \cm^\infty$, then $\frac{f_i}{||x||}\in C^\infty(\R^p,o)$.
 And thus $\frac{f_1  f_2}{||x||}\in (f_1)\cap (f_2)$ but  $\frac{f_1  f_2}{||x||}\not\in (f_1 f_2)$. Thus any two elements of $\cm^\infty$ are not coprime.
\eei

For the ring of continuous functions,  $R=C^0(\cU)$, the theorem is  useless, as  any elements $f_1,f_2\in  R$ with $V(f_1)\cap V(f_2)\neq\empty$
   are not co-prime.   Indeed, both $f_1$ and $f_2$ are divisible by $\sqrt{|f_1|+|f_2|}$, thus
  $\frac{f_1\cdot f_2}{\sqrt{|f_1|+|f_2|}}\in (f_1)\cap(f_2)\sset R$, thus $(f_1)\cap(f_2)\supsetneq (f_1\cdot f_2)$.
   See \cite{Grove-Pedersen} for diagonalization criteria in this case.
\eee

\subsection{Decomposition of determinantal representations/sheaves on plane curves}\label{Sec.Decomposition.Det.Reps.Curves}
Let $ R=\k[x_0,\dots,x_{l}]$,   with $l\ge2$  and  the algebraically closed field, $\k=\bar\k$. Thus $\Proj( R)=\P^l_\k$.
\bee[\bf i.]
\item  Suppose the matrix is homogeneous, $A\in \Mats{R_d}$,
   thus $A$ is a (non-linear) determinantal representation of the projective hypersurface $V(\det(A))\sset \P^l$. Assume $\det(A)=f_1\cdot f_2$, co-prime elements.
    Corollary \ref{Thm.Decompos.Proj.vs.Local.square.case} gives:
        $A$ is decomposable (as an $R$-matrix) \iff  for each (closed) point  $V(\cm)\in \P V(f_1)\cap \P V(f_2)\sset \P^l$
         the localized version $A^{(\cm)}$ is decomposable.
   For $d=1$,    $char(\k)=0$, this recovers  \cite[Theorem 3.1]{Kerner-Vinnikov.Det.Reps.Global}.

\item  For $l=2$ we get a reducible determinantal curve $\P V(det(A))=C_1\cup C_2\sset\P^2$.
  The curves $C_1$, $C_2$ intersect at a finite (non-zero) number of points.
 Assume:
 \bei
 \item  all these points are ordinary multiple points of $C_1\cup C_2$, see \cite[pg.38]{Hartshorne};
 \item  $A^{(\cm)}$ is a locally maximal determinantal representation
   at each such point, see corollary \ref{Thm.Cor.4.7} and example
 \ref{Ex.Decomposability.Square.Matrices.Maximal.Corank.origin}.
 \eei
 Then $A^{(\cm)}$ decomposes locally at each such point.
  Therefore  $A$ is globally $(f_1,f_2)$-decomposable.

\item(Torsion-free sheaves on reducible plane curves) Take a plane curve $C=\cup C_j\sset \P^2_\k$, here  each $C_j$ can be further reducible, but is
 reduced, and the intersections $\{C_j\cap C_i\}_{i\neq j}$ are finite. Take its determinantal representation,
 i.e.,  a homogeneous matrix (see \S\ref{Sec.Background.Notations.Conventions}.v) $A\in \Mats{\k[x_0,x_1,x_2]}$ with $det(A)=\prod f_j$.
   Then $\Coker(A)$ is a coherent sheaf on the curve $C$.
  It is torsion free and  of rank one on each $C_j$.

 Suppose all the intersection points of $\{C_i\cap C_j\}$ are nodes. Then the stalk of $\cF$ at each such point is either a free module or decomposes.
 Therefore we get:
\[
\cF\cong \oplus \cF_j \text{ with } Supp(\cF_j)=C_j\hspace{0.5cm}  \text{ \iff } \hspace{0.5cm} \cF  \text{ is
  not locally free at each such node.}
\]
More generally, suppose for each point $p\in C_i\cap C_j$ the curve-germs have no common tangents, i.e. their tangent cones intersect properly, $T_{(C_i,p)}\cap T_{(C_j,p)}=(0).$ Suppose the determinantal representation is maximal at each point, i.e. $ord(C,p)=dim_\k ker(A|_p)$.
 Then $\cF$ decomposes locally at each point of $\{C_i\cap C_j\}$, and thus decomposes globally.
  \eee

\subsection{Simultaneous diagonal reduction of tuples of matrices/linear determinantal representations}\label{Sec.Decomposition.Tuples.of.Matrices}

 Take an algebraically closed field,  $\k=\bar\k$, and a (possibly infinite) collection of matrices $\{A_\al\in \Mat{\k}\}_\al$.
 The classical question is to ensure the simultaneous diagonalization of this tuple.

    We can assume $\cap \ker(A^\al)=0\sset \k^n$, otherwise one restricts to a complementary subspace of $\cap \ker(A^\al)$ in $\k^n$.

 Split the sizes, $m=\sum m_i$, $n=\sum n_i$. We are looking for $U\in GL(m,\k)$, $V\in GL(n,\k)$ to ensure:
 \beq
 U\cdot A_\al\cdot V^{-1}=\oplus_i A_\al^i ,\quad \{A_\al^i\in \operatorname{Mat}_{m_i\times n_i}(\k)\}_\al,\quad \forall \al.
 \eeq
 First we formulate this as the decomposition problem of one matrix, over a ring.
 Take the vector space $\prod \k\bl x_\al\br$ (we allow infinite linear combinations)
  and the graded ring $R=\k[\prod \k\bl x_\al\br]$.  In the finite case this is just $\k[\ux]$, otherwise $R$ is non-Noetherian.
   Take the matrix $A:=\sum x_\al\cdot A_\al\in \Mat{R}$.
 The collection $\{A_\al\}$ is simultaneously decomposable \iff $A$ is decomposable by
  $GL(m,R)\times GL(n,R)$,
   see lemma \ref{Thm.Decomposability.Graded-to-Local.reduction} . (We can even pass to the localization, $R_{R_{>0}}$.)

 \medskip

In the same way one transforms the simultaneous block-diagonalization (of tuples of matrices) by congruence/conjugation into that for $A\in \Mats{R}$.
  More generally, the decomposability questions of quiver representations are reduced to those of matrices over a ring, see \cite{Kerner-Block-Diag.II}.

Below we work with finite tuples of square matrices.

\subsubsection{Reduction from $\{A_1,\dots,A_{c\ge4}\}$ to triples of matrices} Take linear forms, $l(\uA)=\sum^c_{i=1} a_i A_i.$
 \bcor   Suppose the  generic such linear forms $l_1,l_2,l_3$ satisfy:
\bei
\item   $det[\sum^3_{j=1}y_j l_j(\uA)]=\tf_1(y)\cdots \tf_r(y)$, pairwise coprime polynomials in $y=(y_1,y_2,y_3)$,
\item and  the matrix $\sum^3_{j=1}y_j\cdot l_j(\uA)$ is $(\tf_1,\dots,\tf_r)$-decomposable.
\eei
   Then the tuple $\{A_1,\dots,A_c\}$ is simultaneously decomposable, $\{A_1,\dots,A_c\}\sim\oplus^r_{j=1} \{A^{(j)}_1,\dots,A^{(j)}_c\}.$
 \ecor
\bpr
The passage $\Mats{\k[x]}\ni \sum x_i A_i\rightsquigarrow \sum^3_{j=1}y_j l_j(\uA)=\sum^c_{i=1}A_i l^t_i(y)$ can be realized by the quotient map
 $\k[x]\to\quots{\k[x]}{(l_4(x),\dots,l_c(x))}$. Here $l_4,\dots,l_c$ are linear forms, and the plane $V(l_4,\dots,l_c)\sset \k^c$ is generic for
  the given hypersurface $V(\det[\sum x_i A_i])\sset \k^c$.

This generic plane section is reducible, $\det[\sum x_i A_i]|_{V(l_4,\dots,l_c)}=\tf_1(y)\cdots \tf_r(y)\in \k[y_1,y_2,y_2]$. Therefore by Bertini theorem
 (in any characteristic) \cite[pg.179]{Hartshorne}, we get:
 $\det[\sum x_i A_i]= f_1(x)\cdots  f_r(x)$. Moreover, $ f_1(x),\dots,f_r(x)$ are pairwise coprime, as $\tf_1(y),\dots,\tf_r(y)$ are.

Now we would like to invoke corollary \ref{Thm.Decompos.over.R.vs.R/a.Square.Case}. Thus we should verify: $\bR\cdot (Im A:f_j)\stackrel{?}{=}Im \bA:\bar{f}_j.$   Here $f_j\in R:=\k[x_1\dots x_c]$ and $\bar{f}_j\in \bR:=\quots{R}{(l_4(x),\dots,l_c(x))}.$
\\
Denote $M\!:=\!Coker(A).$ Observe: $0\!\to\! Tor^R_1(\quots{R}{(f_j),M})\!\to\! M\stackrel{f_j\cdot }{\to}M\!\to\! \quots{M}{f_j M}\!\to\!0.$
  Thus $Im A\!:f_j\!=\!Tor^R_1(\quots{R}{(f_j),M}).$ Similarly:  $Im \bA:\bar{f}_j\!=\!Tor^\bR_1(\quots{\bR}{(\bar{f}_j),\bbM}).$
   Thus we should verify: $\bR\otimes Tor^R_1(\quots{R}{(f_j),M})\cong Tor^\bR_1(\quots{\bR}{(\bar{f}_j),\bbM}).$

Applying $\bR\otimes$ to  $0\to Tor^R_1(\quots{R}{(f_j),M})\to M\stackrel{f_j\cdot }{\to}f_j M\to 0$ we get:
 \beq
 Tor^R_1(\bR,f_jM)\to \bR\otimes Tor^R_1(\quots{R}{(f_j),M})\to \bbM\stackrel{\bar{f}_j\cdot }{\to}\bar{f}_j \bbM\to 0
 \eeq
Compare this to  $0\to Tor^\bR_1(\quots{\bR}{(\bar{f}_j),\bbM})\to \bbM \to \bar{f}_j \bbM\to 0.$
 Therefore it is enough to verify: $Tor^R_1(\bR,f_jM)=0.$ Equivalently: the sequence $l_4,\dots,l_c$ is $f_j M$-regular, for generic linear forms $l_4,\dots,l_c.$ Equivalently: $depth_{(x)}f_j M\ge c-3.$

But $depth_{(x)}f_j M=depth_{(x)}M=c-1.$ Hence the statement.
\epr

\subsubsection{The case of a pair of matrices}
\bcor
Let $A_1,A_2\in \Mats{\k}$ with eigenvalues $\{\la^{(1)}_\bullet\}$ and $\{\la^{(2)}_\bullet\}$. Take the
 corresponding generalized eigenspaces    $\{V^{(1)}_{\la^{(1)}_\bullet}\}$, $\{V^{(2)}_{\la^{(2)}_\bullet}\}$.
  If $dim[V^{(1)}_{\la^{(1)}_i}\cap V^{(2)}_{\la^{(2)}_j}]\le 1$
 for any $i,j,$ then the pair $\{A_1,A_2\}$ admits simultaneous diagonal reduction.
\ecor
\bpr
The homogeneous polynomial splits, $\det[x_1 A_1+x_2A_2]=\prod^n_{k=1} l_k(x_1,x_2).$
 We claim: all the forms $l_\bullet$ are  pairwise-coprime. Indeed, we can put e.g. $A_1=-\one,$
  then  $\det[x_1 A_1+x_2A_2]$ is the characteristic polynomial. And by our assumption it has no multiple roots.

Therefore the homogeneous polynomials
$\{\frac{det(x_1 A_1+x_2A_2)}{l_k(x_1,x_2)}\}_k$ are $\k$-linearly independent. Hence
 $I_{n-1}(x_1 A_1+x_2A_2)\sseteq (x_1,x_2)^{n-1}= \sum_k (\frac{det(x_1 A_1+x_2A_2)}{l_k(x_1,x_2)}).$
  Thus $x_1A_1+x_2A_2$ is $\{l_k\}$-decomposable.

\epr

\subsubsection{The case of  matrix-triples}
  For each point  $(x_0:x_1:x_2)\in \P^2_\k$ take the kernel $Ker(\sum x_iA_i)\sseteq \k^n$. This space is non-trivial
   \iff $(x_0:x_1:x_2)\in C:=V(det\sum x_iA_i)\sset \P^2$, the associated determinantal curve.
 Recall that $dim_\k\big(Ker \sum x_iA_i \big)=1$ at the smooth points of $C$.
 Call  $(x_0,x_1,x_2)\in \P^2_\k$ `an essential singular point' if $dim_\k\big( Ker \sum x_iA_i \big)\ge2$.
 An essential singular point is necessarily a singular point  for $C.$  There can be also non-essential singular points of $C,$ called `accidental', see \cite[pg.1622]{Kerner-Vinnikov.Det.Reps.Global}.

   A triple admits the diagonal reduction exactly when the total dimension of such kernels at essential singular points is the maximal possible.
    (And then the curve $C$ is a line-arrangement.)
\bprop
 Assume that the polynomial $\det[\sum x_iA_i]\in \k[x_0,x_1,x_2]$ is square-free.
 \bee[\bf 1.]
 \item Then $\sum_{p\in \P^2} \bin{dim_\k (Ker \sum x_i A_i|_p )}{2}\le \bin{n}{2}.$
 \item
The triple $\{A_0,A_1,A_2\}$ admits the diagonal reduction \iff $\sum_{p\in \P^2} \bin{dim_\k (Ker \sum x_i A_i|_p )}{2}=\bin{n}{2}$.
 \eee
\eprop
Here the sum   (over all the points of $\P^2$) is finite, as only the (essential) singular points $p\in \P^2$ of $C$ contribute.
\bpr (A preparation: the  basic invariant  of curve singularities.) Take the normalization of a
 reduced (possibly reducible) curve germ,
 $\amalg^r_{i=1}(\tC_i,p_i)\to (C,p)$. This defines the embedding of the local rings, $\cO_{(C,p)}\hookrightarrow \prod_i \cO_{(\tC_i,p_i)}$.
  Then the delta invariant is the vector space dimension of the quotient, $\de(C,p):=dim_\k\quots{\prod_i \cO_{(\tC_i,p_i)}}{\cO_{(C,p)}}.$
  See  Exercise 1.8 in \S IV.1 of \cite{Hartshorne}, or pg. 206 of \cite{G.L.S.}

  Recall the basic properties of $\de$.   Let the multiplicity of $(C,p)$ be $m,$ then:
\bei
\item (Proposition 3.34 of \cite{G.L.S.})   $\de(C,x)\ge\bin{m}{2};$ moreover, $\de(C,x)=\bin{m}{2}$ \iff $(C,x)$ is an ordinary multiple point.
\item (Theorem 2.54 of \cite{G.L.S.}) $\de$ is non-increasing in small deformations of the curve-germ.
\item  In particular, for a reduced plane projective curve $C\sset \P^2_\k,$ of degree $n,$ one has: $\de(C)\le \bin{n}{2}.$
 (By the degeneration of $C$ to the arrangement of $n$ lines through one point.)
  And if $C$ is a line arrangement (with arbitrary combinatorics), then $\de(C)= \bin{n}{2}.$
\eei

For a projective (not necessarily planar) curve one has the global normalization, $\tC\!=\!\amalg^r_{j=1} \tC_j\!\to\! C=\cup^r_{j=1} C_j$.
   Take the exact sequence of structure sheaves
   \beq
   0\!\to\! \cO_C\!\to\! \cO_\tC\!\to \!\quots{\cO_\tC}{\cO_C}\!\cong \!\oplus_{p\in Sing(C)}(\quots{\cO_\tC}{\cO_C})_p\!\to\!0.
   \eeq
    Then the total delta invariant can be computed by taking cohomology and applying Riemann-Roch, e.g. exercise 1.9 in \S V.1 of \cite{Hartshorne}. One gets:
   \beq
   \de(C):=\sum_{p\in Sing(C)} \de(C,p)=dim_\k \quots{\cO_\tC}{\cO_C}=r-1+h^1(\cO_C)-h^1(\cO_{\tC}).
   \eeq
Here $h^1(\cO_{\tC})=\sum g(\tC_j)$. For a plane curve, $C\sset \P^2$, of degree $n$, we compute $h^1(\cO_C)$:
\beq
0\!\to\! \cO_{\P^2}(-n)\!\to\! \cO_{\P^2}\!\to\! \cO_C\!\to\! 0,   \text{ hence }   h^1(\cO_C)\!=h^2(\cO_{\P^2}(-n))=h^0(\cO_{\P^2}(n-3))=\bin{n-1}{2}.
\eeq
Altogether, the total delta is: $\de(C)=\bin{n-1}{2}-1-\sum^r_{j=1} (g(\tC_j)-1)$. In particular, $\de(C)\le \bin{n-1}{2}-1+r\le \bin{n}{2}$.
 The equality  holds here \iff the curve is a line arrangement.

\bee[\bf \text{Part} 1.]
\item
For each   point $p  \in \P^2$ the multiplicity of the determinantal curve $V(det( \sum x_iA_i|_p))$ is at least $dim_\k\big( Ker \sum x_iA_i|_p \big)$.
 Therefore the local delta-invariant  of the reduced plane curve singularity
  is $\de(C,p )\ge \bin{dim_\k ( Ker \sum x_iA_i|_{p }  )}{2}$.

Finally, as $C$ is a reduced plane curve, the total $\de$-invariant   is bounded by that of the line arrangement,
  i.e.,  $\de(C)=\sum_p \de(C,p )\le \bin{n}{2}$.
   Therefore $\sum_p \bin{dim_\k (Ker \sum x_i A_i|_p )}{2}\le \bin{n}{2}.$
\item
{\bf The part $\Lleftarrow$.}
The condition $\sum \bin{dim_\k (Ker \sum x_i A_i )}{2}=\bin{n}{2}$ implies $\de(C)=\sum \de(C,p_\al)= \bin{n}{2}$.
 Thus  the curve $V(det\sum x_iA_i)\sset \P^2$ is a (reduced) line arrangement, i.e.,  $\det(\sum x_i A_i )$ splits into pairwise
  independent linear forms. Therefore all the singular points of $C$ are   ordinary multiple points, and for each point
    $\de(C,p_\al)=\bin{dim_\k (Ker \sum x_i A_i |_{p_\al})}{2}.$ Thus
   $Ker( \sum x_i A_i|_{p_\al})$ is of maximal possible dimension, i.e.,  $ \sum x_i A_i|_{p_\al}$ is of maximal corank.
  Then we get the local decomposability, by example \ref{Ex.Decomposability.Square.Matrices.Maximal.Corank.origin}.
  Finally, as $\sum x_iA_i$  decomposes locally at all the intersection points, it
    decomposes globally, by corollary \ref{Thm.Decompos.Proj.vs.Local.square.case}.

\noindent {\bf The part $\Rrightarrow$.} If the triple $\{A_0,A_1,A_2\}$ admits the diagonal reduction then the (square-free) polynomial $\det(\sum x_iA_i)$ splits into pairwise independent linear factors.
 Thus the curve is a (reduced) line arrangement.
 The multiplicity of this curve at a point $p_\al$ is exactly
  $m_\al:=dim_\k\big( Ker \sum x_iA_i |_{p_\al}\big)$, and the local delta-invariant equals $\bin{m_\al}{2}$.  Finally, as the curve is a line arrangement,
   the total delta invariant equals $\sum \bin{m_\al}{2}=\bin{n}{2}$.
\epr\eee

\end{document}